\RequirePackage{fix-cm}
\documentclass[smallextended]{svjour3}   
\smartqed      
\smartqed  
\usepackage{mathptmx}     
\usepackage[utf8]{inputenc}
\usepackage[british]{babel}
\usepackage[T1]{fontenc}
\usepackage{amsfonts}
\usepackage{color}
\usepackage{bm}
\usepackage{amsmath}

\newcommand{\half}{\frac{1}{2}}
\newcommand\diff[2]{\frac{\textrm{d}{#1}}{\textrm{d}{#2}}}
\newcommand{\inv}[1]{\frac{1}{#1}}

\def\R{\mathbb{R}}
\def\N{\mathbb{N}}

\def\bZ{{\mathbf Z}}

\def\bC{{\mathbf C}}

\def\S{\mathbb{S}}
\def\0{{\bf 0}}
\def\x{{\bf x}}

\def\rE{\mathrm{E}}     

\def\S{\mathbb{S}}
\def\0{{\bf 0}}
\def\x{{\bf x}}

\def\rE{\mathrm{E}}

\def\Md{\mathbb{M}^d}

\def\H{\mathbb{H}}
\def\C{\mathbb{C}}

\def\0{{\bf 0}}
\def\rE{\mathrm{E}}     
\def\cov{\mathop{\rm cov}\nolimits}    

\usepackage{amsfonts, amssymb, latexsym}

\begin{document}

\title{Isotropic Covariance Matrix Functions  on Compact Two-Point Homogeneous Spaces}

\author{  Tianshi Lu   \and  Chunsheng Ma     }

\institute{  T. Lu  \at
              Department of Mathematics, Statistics, and Physics, Wichita State University, Wichita, Kansas 67260-0033, USA  \\
              \email{lu@math.wichita.edu }        
           \and
           C. Ma \at
              Department of Mathematics, Statistics, and Physics, Wichita State University, Wichita, Kansas 67260-0033, USA  \\
            \email{chunsheng.ma@wichita.edu}
}

\date{May 5, 2019} 

\maketitle

\begin{abstract}
 The  covariance matrix function is characterized in this paper for a Gaussian or elliptically contoured vector random field that is stationary, isotropic, and mean square continuous on the compact two-point homogeneous space.
 Necessary and sufficient conditions are derived for a symmetric and continuous matrix function to be an isotropic covariance matrix function on all
 compact two-point homogeneous spaces.
 It is also shown that,
 for a symmetric and continuous matrix function with compact support,
 if it makes an isotropic covariance matrix function in the Euclidean space, then it  makes an isotropic covariance matrix function
 on the sphere or the real projective space.



\end{abstract}

\section{Introduction}

A  $d$-dimensional compact two-point homogeneous space  $\mathbb{M}^d$   is a compact Riemannian symmetric space of rank one, and  belongs to one of the
  following five families  (\cite{Helgason2011},   \cite{Wang1952}): the unit spheres $\S^d$ ($ d =1, 2, \ldots$), the real projective spaces $\mathbb{P}^d(\R)$ ($d = 2, 4, \ldots$),
  the complex projective spaces  $\mathbb{P}^d(\mathbb{C})$ ($d = 4, 6, \ldots$), the quaternionic  projective spaces
  $\mathbb{P}^d(\mathbb{H})$ ($d = 8, 12, \ldots$), and the Cayley  elliptic plane $\mathbb{P}^{16}  (Cay)$ or $\mathbb{P}^{16}  (\mathbb{O})$.
There are at least two different approaches to the subject of compact two-point homogeneous spaces   \cite{MaMalyarenko2018},
including an approach based on Lie algebras and a geometric approach, which are used in   probabilistic literature \cite{Askey1976}, \cite{Gangolli1967},  \cite{Malyarenko2013}, statistical literaure  \cite{Patrangenaru2016}, and  approximation theory literature \cite{Azevedo2017}, \cite{BrownDai2005}.
All compact two-point homogeneous spaces share the same property  that all   geodesics in a given one of these spaces are closed and have the same length \cite{Gangolli1967}. In particular, when the unit sphere $\mathbb{S}^d$ is embedded into the space $\mathbb{R}^{d+1}$, the length of any geodesic line is equal to that of the unit circle, that is, $2\pi$. In what follows, the distance $\rho (\x_1, \x_2)$ between two points $\x_1$ and $\x_2$ on $\Md$ is defined in such a way that  the length of any geodesic line on all $\mathbb{M}^d$ is equal to $2\pi$, or 
the distance between any two points is bounded between 0 and $\pi$, {\em i.e.},
$0 \le \rho (\x_1, \x_2) \le  \pi$. 
Over $\S^d$, for instance, $\rho (\x_1, \x_2) $ is defined by $\rho (\x_1, \x_2)= \arccos (\x_1' \x_2), \x_1, \x_2 \in \S^d$, where $\x_1'\x_2$ is the inner product between $\x_1$ and $\x_2$.
Expressions of
$\rho (\x_1, \x_2)$ on other spaces may be found in \cite{Bhattacharya2012}.

   Gaussian random fields on $\Md$ have been studied  in \cite{Askey1976},     \cite{Gangolli1967},   \cite{Malyarenko2013},
  among others,
  while theoretical investigations and
practical applications of  scalar and vector random fields on spheres may be found in
\cite{Askey1976}, \cite{Bingham1973},  \cite{Cheng2016},  \cite{Cohen2012}, \cite{Dovidio2014}
              \cite{Gangolli1967},
             \cite{Leonenko2012},  \cite{Leonenko2013},  \cite{Ma2015}-\cite{Ma2017},
         \cite{Malyarenko2013},  \cite{Malyarenko1992},
    \cite{Yadrenko1983}-\cite{Yaglom1987}.
    Recently, a series representation is presented in \cite{MaMalyarenko2018} for a vector random field that is isotropic and mean square continuous on
    $\Md$ and stationary on a temporal domain, and
a general form of the covariance matrix function is derived  for such a vector random field, which involve Jacobi polynomials and the distance defined on $\Md$.
It is called for parametric and semiparametric covariance matrix structures on $\Md$  in \cite{MaMalyarenko2018}, which are the topics of  this paper.

\begin{table}
 \centering
\caption{ Parameters $\alpha$ and $\beta$ associated with Jacobi polynomials over $\Md$}
  \label{tab:1}
    \begin{tabular}{|l| c |  c|}
    \hline
    $\mathbb{M}^d$ &    $\alpha$ &  $\beta$ \\
    \hline
    $\mathbb{S}^d$, $d=1$, $2$, \dots &  $\frac{d-2}{2}$ &  $\frac{d-2}{2}$  \\
    $\mathbb{P}^d(\mathbb{R})$, $d=2$, $3$, \dots &   $\frac{d-2}{2}$  &   $-\frac{1}{2}$  \\
    $\mathbb{P}^d(\mathbb{C})$, $d=4$, $6$, \dots &  $\frac{d-2}{2}$  & 0  \\
    $\mathbb{P}^d(\mathbb{H})$, $d=8$, $12$, \dots &   $\frac{d-2}{2}$  &  1 \\
    $\mathbb{P}^{16}(Cay)$ &    7   &  3 \\
    \hline
  \end{tabular}
\end{table}

Consider an $m$-variate second-order random field $\{ \bZ (\x), \x \in \Md \}$.  It is called   a stationary (homogeneous) and isotropic random field, if its mean function $\rE \bZ (\x) = ( \rE Z_1(\x), \ldots, \rE Z_m (\x) )'$
does not depend on $\x$, and its covariance matrix function,
    $$ \cov (  \bZ (\x_1),  \bZ ( \x_2) ) = \rE [ ( \bZ (x_1) - \rE \bZ (\x_1)) ( \bZ (x_2) - \rE \bZ (\x_2))' ],    ~~~~~~ \x_1, \x_2 \in \Md, $$
    depends only on the distance $\rho (\x_1, \x_2)$ between $\x_1$ and $\x_2$.
   We denote such a covariance matrix function by $\bC( \rho (\x_1, \x_2)), \x_1, \x_2 \in \Md, $
  and call it an isotropic covariance matrix function on $\Md$.
   An isotropic  random field $\{ \bZ (\x), \x \in \Md \}$ is said to be mean square continuous if, for $k =1, \ldots, m$,
     $$ \rE  | Z_k (\x_1) -Z_k (\x_2) |^2 \to 0,  ~~ \mbox{as} ~~ \rho (\x_1, \x_2 ) \to 0,  ~ \x_1, \x_2 \in \mathbb{M}^d. $$
     It implies the continuity of each entry of the  associated covariance matrix function  in terms of $\rho (\x_1, \x_2)$.

     An $m$-variate isotropic and mean square continuous random field on $\Md$ has a series representation \cite{MaMalyarenko2018}, for $d \ge 2$,
        $$ \bZ ( \x ) = \sum_{n=0}^\infty  \mathbf{B}_n^{\frac{1}{2}} \mathbf{V}_n  P_n^{ (\alpha, \beta) } ( \cos \rho (\vec{x},  \vec{U} )),
      ~~~~~~  \x   \in  \Md, $$
     where
 $\{ \mathbf{V}_n,  n  \in \mathbb{N}_0 \}$ is a sequence of  independent $m$-variate random vectors with
$\rE ( \mathbf{V}_n)= \0$ and  $\cov ( \mathbf{V}_n, \mathbf{V}_n ) =  a_n^2 \mathbf{I}_m$,     $\mathbf{U}$  is a
        random vector uniformly distributed on $\mathbb{M}^d$  and is independent of
        $\{\, \mathbf{V}_n,  n  \in \mathbb{N}_0\, \}$,  $\{ \mathbf{B}_n,  n \in \mathbb{N}_0 \}$ is a sequence of $m \times m$ positive definite matrices,
        $\sum\limits_{n=0}^\infty  \mathbf{B}_n  P_n^{(\alpha, \beta) } \left( 1 \right)$ converges, $\mathbf{I}_m$ is an $m \times m$ identity matrix, $\mathbb{N}_0$ and $\N$ denote the sets of nonnegative integers and of positive integers, respectively,
        \begin{equation}
        \label{JacobiPolynomial}
          P_n^{(\alpha, \beta)} (x) = \frac{\Gamma (\alpha+n+1)}{n! \Gamma (\alpha+\beta+n+1)}\sum_{k=0}^n\binom{n}{k}\frac{\Gamma (\alpha+\beta+n+k+1)}{\Gamma ( \alpha+k+1 )} \left(\frac{x-1}{2} \right)^k,
          \end{equation}
         \hfill $ \quad x \in \R, \quad n \in \mathbb{N}_0, $

\noindent
        are Jacobi polynomials \cite{Szego1975} with specific pairs $\alpha$ and $\beta$ given in Table \ref{tab:1}, and
          \begin{equation}
           \label{a.n.definition}
             a_n=\left(\frac{\Gamma(\beta+1)(2 n +\alpha+\beta+1)\Gamma(n+\alpha+\beta+1)}
{\Gamma(\alpha+\beta+2)\Gamma(n+\beta+1)}\right)^{\frac{1}{2}},  \qquad   n \in \mathbb{N}_0.
           \end{equation}
           The  covariance matrix function of $\{ \bZ(\x),  \x \in \Md \}$ is
          \begin{equation}
     \label{cov.mf1}
       \bC( \rho (\x_1, \x_2))  = \sum_{n=0}^\infty  \mathbf{B}_n  P_n^{(\alpha, \beta) } \left( \cos  \rho (\x_1, \x_2) \right),
      ~~~~~~ \x_1, \x_2 \in \mathbb{M}^d.
      \end{equation}
On the other hand,  there exists an $m$-variate isotropic Gaussian
or elliptically contoured random field on $\Md$ with $\bC( \rho (\x_1, \x_2))$  as its covariance matrix function \cite{MaMalyarenko2018},
if  $\bC( \rho (\x_1, \x_2))$ is  an $m \times m$ symmetric matrix function  of the form (\ref{cov.mf1}).

Given a symmetric matrix function $\bC (\vartheta)$ whose entries are continuous on $[0, \pi]$,
Section 2 presents the characterizations for $\bC( \rho (\x_1, \x_2))$ to be the covariance matrix function of an isotropic elliptically contoured vector random field on $\Md$, in terms of the positive definiteness of a sequence of symmetric matrices.
It is  characterized  in Section 3  for   $\bC( \rho (\x_1, \x_2))$
to be an isotropic covariance matrix function  on all possible $\Md$.
If $\bC (\vartheta)$ makes $\bC (\| \x_1-\x_2 \|)$ an isotropic covariance matrix function in $\R^d$,
 does it make $\bC( \rho (\x_1, \x_2))$ an isotropic covariance matrix function on $\Md$?  A partial answer to this question or the conjecture in
  \cite{NieMa2019}  is given in Section 4,  when
$\bC (\vartheta)$ is compactly supported on $[0, \pi]$ and $d$ is odd.
Proofs of theorems are given in Section 5.

 \section{Isotropic covariance  matrix functions  on $\Md$  }

 The covariance matrix function is characterized in this section of an isotropic and mean square continuous elliptically contoured vector random field on $\Md$.     Theorem \ref{thm1}
  provides a useful tool for verifying whether a continuous matrix function is the covariance matrix function of an isotropic
  vector elliptically random field   on $\Md$, by checking that each of a sequence of matrices is positive definite and a relevant infinite series is convergent, and  Theorem \ref{thm2} presents  the interrelationship of an isotropic covariance matrix function on different compact two-point homogeneous spaces.

   \begin{theorem}
   \label{thm1}
    Let  $\alpha$ and $\beta$ be the pair for $\Md$ in Table \ref{tab:1}.
    For an $m \times m$ symmetric matrix function $\bC (\vartheta)$ whose entries are continuous on $[0, \pi]$,
    the following statements are equivalent:

     \begin{itemize}
     \item[\textup{(i)}]
    $\bC (\rho (\x_1, \x_2))$ is the covariance matrix function of an $m$-variate isotropic elliptically contoured random field on $\Md$;

     \item[\textup{(ii)}] $\bC (\vartheta)$ is of the form
        \begin{equation}
        \label{cov.mf2}
        \bC (\vartheta) = \sum_{n=0}^\infty \mathbf{B}_n P_n^{ (\alpha, \beta)} ( \cos \vartheta),  ~~~~~ \vartheta \in [0, \pi ],
        \end{equation}
        where $\{ \mathbf{B}_n, n \in \N_0 \}$ is a sequence of $m \times m$ positive definite matrices, and the series
         $\sum\limits_{n=0}^\infty  n^{\alpha} \mathbf{B}_n  $ converges;

      \item[\textup{(iii)}]  the matrices
      \begin{equation}
      \label{thm1.eq1}
      \mathbf{H}_n^{(\alpha, \beta) } = \int_0^\pi \bC (\vartheta)  P_n^{(\alpha, \beta) } \left( \cos  \vartheta \right)
                 \sin^{2 \alpha+1}  \left( \frac{\vartheta}{2} \right) \cos^{2 \beta+1}  \left(  \frac{\vartheta}{2}  \right) d \vartheta,
             ~~~~ n \in \mathbb{N}_0,
      \end{equation}
      are positive definite, and the series $\sum\limits_{n=0}^\infty  n^{\alpha+1} \mathbf{H}_n^{(\alpha, \beta) }  $ converges.
      \end{itemize}
   \end{theorem}

Note that $P_n^{(\alpha, \beta)} (1) = \frac{\Gamma (n+\alpha+1)}{\Gamma (n+1) \Gamma (\alpha+1)}$, $n \in \N_0$.
  By the asymptotic formula (5.11.12) of \cite{Olver2010},
           $ \frac{\Gamma (n +\alpha+1)}{\Gamma (n+1)} \sim n^\alpha$ ~ ($n \to \infty$),
             $\sum\limits_{n=0}^\infty \mathbf{B}_n P_n^{(\alpha, \beta)} (1)$ converges if and only if
           $\sum\limits_{n=0}^\infty n^\alpha \, \mathbf{B}_n$ converges.
The convergence of $\sum\limits_{n=0}^\infty n^\alpha \mathbf{B}_n$ in Theorem \ref{thm1} (ii)  is necessary to guarantee the convergence of the series in   (\ref{cov.mf2}) for $\vartheta=0$, and is also sufficient for the convergence for all $\vartheta \in [0, \pi]$, since $|P^{(\alpha,\beta)}_n(\cos\theta)|\le P^{(\alpha,\beta)}_n(1),  n \in \N_0$. The condition that
$\sum\limits_{n=0}^\infty  n^{\alpha+1} \mathbf{H}_n^{ (\alpha, \beta)}$ converges in Theorem \ref{thm1} (iii)   is equivalent to the convergence of $\sum\limits_{n=0}^\infty n^\alpha \mathbf{B}_n$  in Theorem \ref{thm1} (ii).

    There are two key parameters associated with $\Md$ in Table \ref{tab:1}, $\alpha$ and $\beta$,
    which are not dependent each other, except for $\S^d$ where $\alpha = \beta$.
     The parameter $\beta$  is a constant  with respect to $d$ or $\Md$, except for $\S^d$.
   The following formula expresses   a coefficient  $  \mathbf{H}_n^{(\alpha, \beta) }$ on $\mathbb{M}^d$ in terms of  two coefficients  $ \mathbf{H}_n^{(\alpha-1, \beta) }$ and $ \mathbf{H}_{n+1}^{(\alpha-1, \beta) }$ on $\mathbb{M}^{d-2}$ ($d \ge 3$).

    \begin{corollary}  For $d \ge 3$,
    \begin{equation}
     \label{H.n.idenity}
      \mathbf{H}_n^{(\alpha, \beta) }  = \frac{ (n+\alpha)  \mathbf{H}_n^{(\alpha-1, \beta) }  - (n+1)  \mathbf{H}_{n+1}^{(\alpha-1, \beta) } }{ 2n +\alpha+\beta+1}, ~~~~ n \in \N_0.
     \end{equation}
    \end{corollary}

    Identity (\ref{H.n.idenity}) follows directly from (\ref{thm1.eq1}) and  (4.5.4) of \cite{Szego1975},
       \begin{equation}
       \label{Szego4.5.4}
        \frac{ 2 n+\alpha +\beta +1 }{2}   (1-x) P_n^{(\alpha, \beta)} (x)
        = (n+\alpha)  P_n^{(\alpha-1, \beta)} (x)  - (n+1) P_{n+1}^{(\alpha-1, \beta)} (x), ~~ x \in \R.
        \end{equation}
        A dual identity of (\ref{Szego4.5.4})  is
        \begin{equation}
       \label{Szego4.5.4.dual}
          \frac{ 2 n+\alpha +\beta +1 }{2}   (1+x) P_n^{(\alpha, \beta)} (x)
        = (n+\beta)  P_n^{(\alpha, \beta-1)} (x)  + (n+1) P_{n+1}^{(\alpha, \beta-1)} (x), ~~ x \in \R,
        \end{equation}
        from which and from (\ref{thm1.eq1}) we obtain the following corollary, which is useful over $\S^d$ since $\beta$ keeps  fixed over other spaces.

        \begin{corollary}  For $d \ge 3$, 
    \begin{equation}
     \label{H.n.idenity2}
      \mathbf{H}_n^{(\alpha, \beta) }  = \frac{ (n+\beta)  \mathbf{H}_n^{(\alpha, \beta-1) }  + (n+1)  \mathbf{H}_{n+1}^{(\alpha, \beta-1) } }{ 2n +\alpha+\beta+1},  ~~~~ n \in \N_0.
     \end{equation}
    \end{corollary}

      Notice that the parameter $\alpha$ in Table \ref{tab:1} is either a nonnegative integer or an integer plus half, according to whether the dimension $d$ is even or odd.
For these two cases, in the next two corollaries  we are going to write $  \mathbf{H}_n^{(\alpha, \beta) }$ as a linear combination of  $  \mathbf{H}_j^{(0, \beta) }$ or $  \mathbf{H}_j^{ \left( -\frac{1}{2}, \beta \right) }$,  $j \ge n$, respectively, which are coefficients in low dimensions.

         For an even $d \ge 4$ or a positive integer $\alpha =\frac{d-2}{2}$,
     successively using identity (\ref{Szego4.5.4}), the $n+\alpha$   degree
         polynomial $(1-x)^\alpha  P_n^{(\alpha, \beta)} (x)$ can be expressed as a linear combination of polynomials
        $  P_{n+j}^{(0, \beta)} (x)$, $ j = 0, 1, \ldots,  n+\alpha$. More precisely,     it can be established by induction on $\alpha$ that
        \begin{equation}
        \label{Jacobi.identity1}
        \left(\frac{1-x}{2}\right)^\alpha  P_n^{(\alpha, \beta)} (x) = \sum_{j=0}^\alpha (-1)^j a_j^{(0)}  (n) P_{n+j}^{(0, \beta)} (x), ~~~ x \in \R, ~ n \in \N_0,
        \end{equation}
        where
          \begin{equation}
        \label{Jacobi.identity1.coeff}
        \begin{array}{llr}
       a_j^{(0)}  (n)    & = &    \frac{\alpha! \Gamma(n+\alpha+1) (2n+2j+\beta+1) \Gamma(2n+j+\beta+1)}{
                                  j! (\alpha-j)! n!  \Gamma(2n+j+\alpha+\beta+2) },    \\
                   &   &
                                   ~~~ j =0, 1, \ldots, \alpha.
         \end{array}
        \end{equation}
        The following corollary is derived from (\ref{thm1.eq1}) and (\ref{Jacobi.identity1}).

         \begin{corollary}  For $\alpha \in \N$,
    \begin{equation}
     \label{H.n.idenity3}
      \mathbf{H}_n^{(\alpha, \beta) }  =  \sum_{j=0}^\alpha (-1)^j   a_j^{(0)}  (n)   \mathbf{H}_{n+j}^{(0, \beta) }, ~~~~ n \in \N_0.
     \end{equation}
    \end{corollary}

       For an odd  $d \ge 3$,  $\alpha +\frac{1}{2} =\frac{d-1}{2}$ is a positive integer.
     Successively using identity (\ref{Szego4.5.4}), the $n+\alpha +\frac{1}{2}$   degree
         polynomial $(1-x)^{\alpha+\frac{1}{2}}  P_n^{(\alpha, \beta)} (x)$ can be expressed as a linear combination of polynomials
        $  P_{n+j}^{ \left( -\frac{1}{2}, \beta \right)} (x)$, $ j = 0, 1, \ldots,  n+\alpha+\frac{1}{2}$,     and,  by induction on $\alpha+\frac{1}{2}$,
        \begin{equation}
        \label{Jacobi.identity2}
        \left(\frac{1-x}{2}\right)^{ \alpha+\frac{1}{2} }   P_n^{(\alpha, \beta)} (x) = \sum_{j=0}^{\alpha+\frac{1}{2}}  (-1)^j  a_j^{ \left(  - \frac{1}{2} \right)} (n) P_{n+j}^{ \left( -\frac{1}{2}, \beta \right)} (x), ~~~ x \in \R, ~ n \in \N_0,
        \end{equation}
        where
          \begin{equation}
        \label{Jacobi.identity2.coeff}
        \begin{array}{llr}
       a_j^{ \left(  - \frac{1}{2} \right)} (n)    & = &   \frac{\Gamma\left(\alpha+\frac{3}{2}\right) (n+j)! \Gamma(n+\alpha+1)
       \left(2n+2j+\beta+\frac{1}{2}\right) \Gamma\left(2n+j+\beta+\frac{1}{2}\right) }{
        j! \Gamma\left(\alpha-j+\frac{3}{2}\right) n! \Gamma\left(n+j+\frac{1}{2}\right) \Gamma(2n+j+\alpha+\beta+2) }, \\
                   &   &
                                   ~~~ j =0, 1, \ldots, \alpha.
         \end{array}
        \end{equation}
        The following corollary follows directly from (\ref{thm1.eq1}) and (\ref{Jacobi.identity2}).

         \begin{corollary}  For $\alpha +\frac{1}{2} \in \N$,
    \begin{equation}
     \label{H.n.idenity4}
      \mathbf{H}_n^{(\alpha, \beta) }  =  \sum_{j=0}^{\alpha +\frac{1}{2}}  (-1)^j  a_j^{ \left(  - \frac{1}{2} \right)} (n)    \mathbf{H}_{n+j}^{ \left( -\frac{1}{2}, \beta \right) }, ~~~~ n \in \N_0.
     \end{equation}
    \end{corollary}

     Second-order elliptically contoured  random fields  form one of the largest sets, if not the largest set,
          which allows any possible correlation structure \cite{Ma2011}.
          Examples of elliptically contoured random fields include  Gaussian, Student's t,  Cauchy,  Laplace,
          logistic,  hyperbolic, hyperbolic secant,   variance Gamma, normal inverse Gaussian,   K-differenced,  stable, Linnik,  and Mittag-Leffler random fields.
          The characterizations in Theorem \ref{thm1} are available
          for a  second-order elliptically contoured vector random field.
          However, they may not be available for other non-Gaussian  random fields, such as  a log-Gaussian, $\chi^2$, binomial-$\chi^2$, K-distributed, or skew-Gaussian one,  for which  admissible correlation structure must be investigated on a case-by-case basis.

          In what follows, every  covariance matrix function is set up under the elliptically contoured background.
          To distinguish the distances of  the five families listed in Table 1, whenever necessary, we adopt the symbol  $\rho_{\tiny{\S^d}}( \x_1, \x_2)$ for the distance over $\S^d$,  $\rho_{\tiny{\mathbb{P}^d( \R)}}( \x_1, \x_2)$ for the distance on $\mathbb{P}^d( \R)$, and so on.  
          The next theorem shows the interrelationship of an isotropic covariance matrix function on different compact two-point homogeneous spaces.
           It looks like that   isotropic covaraince matrix structures on $\S^d$ are  richer than those on other    compact two-point homogeneous spaces.

     \begin{theorem}
   \label{thm2}
    Suppose that   $\bC (\vartheta)$ is  an $m \times m$ symmetric matrix function and  each of its entries is continuous on $[0, \pi]$.
      \begin{itemize}
     \item[\textup{(i)}] For an odd $d \ge 3$,
     if $\bC (\vartheta)$ makes
    $\bC \left( \rho_{\tiny{\mathbb{P}^d( \R)}}( \x_1, \x_2) \right)$  an isotropic covariance matrix function on $\mathbb{P}^d (\R)$,
    then  it makes $\bC \left( \rho_{\tiny{\S^d}}( \x_1, \x_2)  \right)$ is an isotropic covariance matrix function on $\S^d$.
    For an even $d$, if    $\bC \left(   \rho_{\tiny{\mathbb{P}^d( \R)}}( \x_1, \x_2)  \right)$ 
      is an isotropic covariance matrix function on $\mathbb{P}^d (\R)$,
    then   $\bC \left(  \rho_{\tiny{\S^{d-1}}}( \x_1, \x_2)  \right)$ is  an isotropic covariance matrix function on $\S^{d-1}$.

      \item[\textup{(ii)}] For an even $d \ge 4$,
     if
    $\bC \left( \rho_{\tiny{\mathbb{P}^d( \C)}}( \x_1, \x_2) \right) $ is an isotropic covariance matrix function  on $\mathbb{P}^d (\C)$,
    then  $\bC \left( \rho_{\tiny{\S^d}}( \x_1, \x_2)  \right)$  is  an isotropic covariance matrix function on $\S^d$ for $ d \ge 4$,  and 
    $\bC \left( \rho_{\tiny{\mathbb{P}^d( \H)}}( \x_1, \x_2) \right)$   is  an isotropic covariance matrix function on $\mathbb{P}^d (\H)$  if $d =8, 12, \ldots.$

    \item[\textup{(iii)}] For  $d = 8, 12, \ldots$,
     if   $\bC \left( \rho_{\tiny{\mathbb{P}^d( \H)}}( \x_1, \x_2) \right)$
     is  an isotropic   covariance matrix function on $\mathbb{P}^d (\H)$,
    then   $\bC \left( \rho_{\tiny{\S^d}}( \x_1, \x_2)  \right)$  is  an isotropic covariance matrix function on $\S^d$.

     \item[\textup{(iv)}] ~  If   $\bC \left( \rho_{\tiny{\S^d}}( \x_1, \x_2)  \right)$  is an isotropic covariance matrix function on $\S^d$, then
     both $\bC \left( \frac{ \rho_{\tiny{\mathbb{P}^d( \R)}}( \x_1, \x_2)  }{2} \right)$  $+\bC \left( \pi - \frac{ \rho_{\tiny{\mathbb{P}^d( \R)}}( \x_1, \x_2)}{2} \right)$
     and 
     
     \noindent  $ \left\{ \bC \left( \frac{\rho_{\tiny{\mathbb{P}^d( \R)}}( \x_1, \x_2)}{2} \right)-\bC \left( \pi - \frac{\rho_{\tiny{\mathbb{P}^d( \R)}}( \x_1, \x_2)}{2} \right) \right\} \cos \left( \frac{\rho_{\tiny{\mathbb{P}^d( \R)}}( \x_1, \x_2) }{2} \right)$
       are  isotropic covariance matrix functions  on $\mathbb{P}^d (\R)$.

     \end{itemize}
     \end{theorem}

   \section{Isotropic covariance matrix functions on all dimensions}

   Except for $\mathbb{P}^{16} (Cay)$, the dimension $d$ of $\Md$ can take infinitely many values, as shown in Table \ref{tab:1}.
   If an $m \times m$  continuous matrix function $\bC (\vartheta)$ on $[0, \pi]$  makes $\bC (\rho (\x_1, \x_2))$
   an isotropic covariance  matrix function on all possible $\Md$ (all five families described in Section 1), then it is called an isotropic covariance matrix function on $\mathbb{M}^\infty$. Such a matrix function is characterized in the following theorem.

   \begin{theorem}
   \label{thm3}
   For an $m \times m$ symmetric matrix function $\bC(\vartheta)$ whose all entries are continuous on $[0, \pi]$,
   the following statements are equivalent:

   \begin{itemize}
   \item[\textup{(i)}] $\bC (\rho (\x_1, \x_2))$ is an isotropic covariance matrix function on $\mathbb{M}^\infty$;

   \item[\textup{(ii)}] $\bC(\vartheta)$ is of the form
    \begin{equation}
   \label{thm3.eq1}
   \bC (\vartheta) = \sum_{n=0}^\infty \mathbf{B}_n (1+ \cos \vartheta )^n, ~~~~~  \vartheta \in [0, \pi],
   \end{equation}
   where  $\{ \mathbf{B}_n, n \in \N_0 \}$ is a  sequence of $m \times m$ positive definite matrices and $\sum\limits_{n=0}^\infty 2^n \mathbf{B}_n$ converges;

    \item[\textup{(iii)}] $\bC \left( \frac{\pi}{2} - \arcsin x \right)$ is of the form
   \begin{equation}
   \label{thm3.eq2}
   \bC \left( \frac{\pi}{2} - \arcsin x  \right) = \sum_{n=0}^\infty \mathbf{B}_n (1+x)^n, ~~~~~ x \in [-1, 1],
   \end{equation}
   where  $\{ \mathbf{B}_n, n \in \N_0 \}$ is a  sequence of $m \times m$ positive definite matrices and $\sum\limits_{n=0}^\infty 2^n \mathbf{B}_n$ converges;

    \item[\textup{(iv)}] $\bC \left( \pi -2 \arcsin x \right)$ is of the form
   \begin{equation}
   \label{thm3.eq3}
   \bC \left( \pi -2 \arcsin x  \right) = \sum_{n=0}^\infty 2^n \mathbf{B}_n x^{2n}, ~~~~~ x \in [0, 1],
   \end{equation}
   where  $\{ \mathbf{B}_n, n \in \N_0 \}$ is a  sequence of $m \times m$ positive definite matrices and $\sum\limits_{n=0}^\infty  2^n \mathbf{B}_n$ converges.
   \end{itemize}
   \end{theorem}

   It is understandable that the characterizations in Theorem \ref{thm3} differ from those on all spheres $\S^\infty$ presented in \cite{Ma2015}.
   Actually, the set of  isotropic covariance matrix functions on $\mathbb{M}^\infty$  is a proper subset of that  on  $\S^\infty$.

    \begin{corollary}
    \label{thm3.cor}
    If $\bC (\rho (\x_1, \x_2))$ is an isotropic covariance matrix function on $\mathbb{M}^\infty$, then

      \begin{itemize}
      \item[\textup{(i)}] $\bC (\vartheta)$ is a positive definite matrix for each fixed $\vartheta \in [ 0, \pi]$;

       \item[\textup{(ii)}] $\bC (\vartheta_1) - \bC (\vartheta_2) $ is a positive definite matrix for  $0 \le \vartheta_1 \le  \vartheta_2 \le  \pi$;

        \item[\textup{(iii)}]  for $\vartheta \in (0, \pi),$ $\bC' (\vartheta)  $ is a negative definite matrix, whenever the derivative exists.
      \end{itemize}

    \end{corollary}

  Since each $\mathbf{B}_n$ is positive definite,    Corollary \ref{thm3.cor}   is due to (\ref{thm3.eq1}),  Part (i) from
  the fact that $1+ \cos \vartheta$ is nonnegative,  Part (ii) from that $1+\cos \vartheta$ is decreasing on $[0, \pi]$, and Part (iii) from
  Part (ii) and
    $\bC' (\vartheta) = - \lim\limits_{\delta \to 0+} \frac{\bC (\vartheta) - \bC (\vartheta+\delta)}{\delta}$.

    \begin{example}
  For  an $m \times m$ symmetric matrix function whose entries are  second order polynomials,
        $$  \bC(\vartheta) = {\bf B}_0 + {\bf B}_1 \vartheta +{\bf B}_2 \vartheta^2,  ~~~~~~~~~~ \vartheta  \in [0, \pi],  $$
      it makes $\bC ( \rho (\x_1, \x_2)) $  an isotropic covariance matrix function on $\mathbb{M}^\infty$ if and only if
                   ${\bf B}_0- \pi^2 {\bf B}_2$ and   $ {\bf B}_2$
                   are positive definite matrices, and ${\bf B}_1= -2 {\bf B}_2 \pi$.
 To derive a form of (\ref{thm3.eq3}) for $\bC  \left( \pi- 2\arcsin x \right) $, we  employ  the Taylor  expansions of $\arcsin x$ and $(\arcsin x)^2$,
  \begin{equation}
         \label{Ex.arcsin}
          \arcsin x = \sum\limits_{n=0}^\infty \frac{(2n)!}{2^{2n} (n!)^2 (2n+1)} x^{2n+1},  ~~~~~~~~~ x \in [-1, 1],   \end{equation}
  \begin{equation}
       \label{Ex.arcsin2}
        (\arcsin x)^2= \sum_{n=1}^\infty  \frac{2^{2n-1} ( (n-1)!)^2}{(2n)!}  x^{2n}, ~~~~~~~~~~ x \in [-1, 1],
     \end{equation}
      and obtain
         \begin{eqnarray*}
         &   & \bC  \left( \pi -2 \arcsin x \right)  \\
                                 & = &    {\bf B}_0 +    \pi   {\bf B}_1     +  \pi^2  {\bf B}_2
                                            - 2 \left( {\bf B}_1+ 2\pi {\bf B}_2 \right) \arcsin x
                                             + 4 {\bf B}_2  (\arcsin  x  )^2 \\
                                  & = &     {\bf B}_0 +    \pi   {\bf B}_1     +  \pi^2  {\bf B}_2
                                            - 2 \left( {\bf B}_1+2 \pi {\bf B}_2 \right)  \sum\limits_{n=0}^\infty \frac{(2n)!}{2^{2n} (n!)^2 (2n+1)}  x^{2n+1}  \\
                                  &  &
                                             + 4 {\bf B}_2  \sum_{n=1}^\infty  \frac{2^{2n-1} ( (n-1)!)^2}{(2n)!}   x^{2n},   ~~ x  \in [0, 1].
         \end{eqnarray*}
         By Theorem \ref{thm3} (iv), $ {\bf B}_0 +    \pi   {\bf B}_1     +  \pi^2  {\bf B}_2 $ and ${\bf B}_2$ must be positive definite, and ${\bf B}_1+2 \pi {\bf B}_2 = \mathbf{0}$.

                Moreover, if    $ {\bf B}_0-\pi^2 {\bf B}_2$  and ${\bf B}_2$ are $m \times m$ positive definite matrices, then
                  $$  \bC( \rho (\x_1, \x_2)) = \left(  {\bf B}_0 - 2 \pi {\bf B}_2 \rho (\x_1, \x_2) +  {\bf B}_2 (\rho (\x_1, \x_2) )^2 \right)^{\circ \ell} $$
                  is  an isotropic covariance matrix function on $\mathbb{M}^\infty$,
          where $\ell$ is a natural number, and ${\bf B}^{\circ \ell}$ denotes the  Hadamard  $\ell$ power of ${\bf B} =(b_{ij})$, whose entries are $b_{ij}^\ell$, the $\ell$ power of $b_{ij}, i, j = 1, \ldots, m$.
          \end{example}

          \begin{example}
          Given two  $m \times m$ symmetric matrices $\mathbf{B}_1$ and  $\mathbf{B}_2$,  consider an $m \times m$ matrix function
   $$  \bC (\vartheta) = \mathbf{B}_1 \exp \left( \frac{\vartheta}{2} \right)+ \mathbf{B}_2 \exp \left( -\frac{\vartheta}{2} \right),    ~~~~~~~ \vartheta \in [0, \pi]. $$
   By Theorem \ref{thm3}, $\bC (\rho (\x_1, \x_2))$ is an isotropic covariance matrix function on $\mathbb{M}^\infty$
    if and only if $\mathbf{B}_2= \mathbf{B}_1 e^\pi$ is a positive definite matrix.
 To see this,  notice that $\exp ( \arcsin  x)$  possesses the Taylor series  with positive coefficients  (see, for instance,  formula 1.216 of \cite{Gradshteyn2007})
      \begin{equation}
      \label{Taylor.series. exp.arcsin}
       \exp ( \arcsin  x) =  \sum_{n=0}^\infty a_n x^n = 1+x+\frac{x^2}{2!}+\frac{2x^3}{3!}+\frac{5x^4}{4!}+\cdots,        ~~~~~~ x \in [-1, 1],
       \end{equation}
    from which we obtain
      \begin{eqnarray*}
      \bC      \left( \pi- 2 \arcsin x \right)
                    &  = &    \sum_{n=0}^\infty  \left(  (-1)^n \mathbf{B}_1 e^{\frac{\pi}{2}}  +\mathbf{B}_2 e^{-\frac{\pi}{2}} \right)  a_n  x^n,  ~~~~~~~~~~ x \in [0, 1].
      \end{eqnarray*}
      A comparison between the last equation with (\ref{thm3.eq3}) results in
      the positive definiteness of $   \mathbf{B}_1 e^{\frac{\pi}{2}}  +\mathbf{B}_2 e^{-\frac{\pi}{2}}$ and
      $$  (-1)^{2 n+1} \mathbf{B}_1 e^{\frac{\pi}{2}}  +\mathbf{B}_2 e^{-\frac{\pi}{2}}  = \mathbf{0},  ~~~~~ n \in \N_0, $$
      or, equivalently,  the positive definiteness of $\mathbf{B}_2= \mathbf{B}_1 e^\pi$.
          \end{example}

             \begin{example}
          Given an $m \times m$ symmetric matrix $\mathbf{B}$ with entries $b_{ij}$, the entries of an $m \times m$ matrix function $\bC(\vartheta)$ are
          defined by
   $$ C_{ij} (\vartheta) = \exp \left(  b_{ij}  \cos \frac{ \vartheta}{2}  \right) + \exp \left( - b_{ij}  \cos \frac{ \vartheta}{2}  \right),    ~~~~~~~ \vartheta \in [0, \pi],  ~~~  i,  j =1, \ldots, m. $$
    It makes $\bC (\rho (\x_1, \x_2))$  an isotropic covariance matrix function on $\mathbb{M}^\infty$
    if and only if ${\bf B}^{\circ 2}$ is a positive definite matrix,  by Theorem \ref{thm3},  since
      $$  \bC      \left( \pi- 2 \arcsin x \right)  = 2 \sum_{n=0}^\infty \frac{ \mathbf{B}^{\circ 2n} }{(2n)!} x^{2n},   ~~~~~ x \in [0, 1], $$
    and  $\sum\limits_{n=0}^\infty \frac{ \mathbf{B}^{\circ 2n} }{(2n)!}$ converges.

     \end{example}

          In the scalar case $m=1$, the following corollary is a consequence of Theorem \ref{thm3},
          and, by Theorem 2 of \cite{Askey1976}, (\ref{thm3.eq4}) below  is an isotropic covariance function on $\mathbb{M}^\infty$.

    \begin{corollary}
   \label{thm3.cor2}
   For a continuous function $C(\vartheta)$  on $[0, \pi]$,
   the following statements are equivalent:

   \begin{itemize}
   \item[\textup{(i)}] $C (\rho (\x_1, \x_2))$ is an isotropic covariance function on $\mathbb{M}^\infty$;

   \item[\textup{(ii)}] $C(\vartheta)$ is of the form
    \begin{equation}
   \label{thm3.eq4}
    C (\vartheta) = \sum_{n=0}^\infty  b_n (1+ \cos \vartheta )^n, ~~~~~  \vartheta \in [0, \pi],
   \end{equation}
   where  $\{ b_n, n \in \N_0 \}$ is a  sequence of  nonnegative numbers and $\sum\limits_{n=0}^\infty 2^n b_n$ converges;

    \item[\textup{(iii)}] $C \left( \frac{\pi}{2} - \arcsin x \right)$ is of the form
   \begin{equation}
   \label{thm3.eq5}
     C \left( \frac{\pi}{2} - \arcsin x  \right) = \sum_{n=0}^\infty b_n (1+x)^n, ~~~~~ x \in [-1, 1],
   \end{equation}
   where  $\{  b_n, n \in \N_0 \}$ is a  sequence of nonnegative numbers and $\sum\limits_{n=0}^\infty 2^n  b_n$ converges;

    \item[\textup{(iv)}] $C \left( \pi -2 \arcsin x \right)$ is of the form
   \begin{equation}
   \label{thm3.eq6}
    C  \left( \pi -2 \arcsin x  \right) = \sum_{n=0}^\infty 2^n b_n x^{2n}, ~~~~~ x \in [0, 1],
   \end{equation}
   where  $\{ b_n, n \in \N_0 \}$ is a  sequence of nonnegative  numbers and $\sum\limits_{n=0}^\infty  2^n b_n$ converges.
   \end{itemize}
   \end{corollary}

  The exponential function $\exp \left(- \frac{\vartheta}{2} \right)$ is  an important generator on all spheres $\S^\infty$. But,  interestingly,  it does not make
  $\exp \left( - \frac{\rho (\x_1, \x_2)}{2} \right)$ an isotropic covariance function on $\mathbb{M}^\infty$, since it follows from  (\ref{Taylor.series. exp.arcsin}) that
       $$ \exp \left(- \frac{\pi - 2 \arcsin x}{2}  \right) = \exp \left(- \frac{\pi}{2} \right) \sum_{n=0}^\infty a_n x^n,   ~~~~~~ 0 \le x  \le  1, $$
       so that (\ref{thm3.eq6}) fails. Nevertheless, $\cosh  \left(  \frac{\pi- \rho (\x_1, \x_2)}{2} \right)$ is an isotropic covariance function on $\mathbb{M}^\infty$, as is seen from Example 2.

   \begin{example}  For a constant $\nu \in (0, 2]$,
      $$ C( \vartheta ) = 1- \left( \sin \frac{\vartheta}{2} \right)^\nu,  ~~~~ \vartheta \in [0, \pi], $$
      makes $C(\rho (\x_1, \x_2))$ an isotropic covariance function on $\mathbb{M}^\infty$.  Corollary \ref{thm3.cor2} is applicable, with a form (\ref{thm3.eq6}) of $C( \pi- 2 \arcsin x )$ given by
        \begin{eqnarray*}
           C  \left( \pi -2 \arcsin x  \right)  & = & 1- \left( \sin \frac{\pi-2 \arcsin x}{2} \right)^\nu \\
                                                              & = & 1- (1-x^2)^{\frac{\nu}{2} } \\
                                                              & = &  \left\{
                                                                  \begin{array}{lll}
                                                                  x^2, ~ &  ~ \nu =2,   \\
                                                                  \sum\limits_{n=1}^\infty \frac{ \prod\limits_{k=1}^n \left( k- \frac{\nu}{2} \right) }{n!} x^{2n},  ~ &  ~ \nu \in (0, 2),  ~ x \in [0, 1).
                                                                  \end{array}  \right.
         \end{eqnarray*}
         Moreover, for constants $\nu_i \in (0, 2]$,
         an $m \times m $ matrix function $\bC (\vartheta)$ with entries
          $$ C_{ij} (\vartheta) = 1- \left( \sin \frac{\vartheta}{2} \right)^{\max (\nu_i, \nu_j)},  ~~~~ \vartheta \in [0, \pi], ~ i, j = 1, \ldots, m, $$
      makes $C(\rho (\x_1, \x_2))$ an isotropic covariance function on $\mathbb{M}^\infty$, by Theorem \ref{thm3},  since an $m \times m $ matrix with entries  \, $k- \frac{\max ( \nu_i, \nu_j)}{2} = \min \left( k- \frac{\nu_i}{2}, k - \frac{\nu_j}{2} \right)$ is positive definite, for $k \ge 1$.

   \end{example}

   \begin{lemma}
   \begin{itemize}
   \item[\textup{(i)}] For $n \in \N_0$,
     \begin{equation}
     \label{Askey1974}
     \left( \frac{1+x}{2} \right)^n =
      \sum_{k=0}^n \varphi_k^{(\alpha,\beta)} P_n^{(\alpha, \beta)} (x),  ~~~~~~ x \in \R,
      \end{equation}
      where
       $$  \varphi_k^{(\alpha,\beta)} =\frac{\Gamma (n+\beta+1)n! (2k +\alpha+\beta+1) \Gamma (k+\alpha+\beta+1) \Gamma (k+\alpha+1)}{ \Gamma (k+n+\alpha+\beta+2) \Gamma (k+\beta+1) k! (n-k)! \Gamma (\alpha+1) P_n^{(\alpha, \beta)} (1)}.
        $$

     \item[\textup{(ii)}]
      With $a_k$ given by (\ref{a.n.definition}),  if $\mathbf{U}$ is a random vector uniformly distributed on $\Md$, then
          $$ Z ( \x) = \sum_{k=0}^n  a_k \left( \varphi_k^{(\alpha,\beta)} \right)^{\frac{1}{2}}  P_n^{(\alpha, \beta)} (  \cos \rho (\x,  \mathbf{U} )  ),
              ~~~~~~ \x \in \Md, $$
          is a scalar isotropic random field on $\Md$ with mean 0 and covariance function $\left( \frac{1+\cos \rho (\x_1, \x_2) }{2} \right)^n$.

        \end{itemize}

   \end{lemma}

       \begin{lemma}
   \begin{itemize}
   \item[\textup{(i)}] For a fixed $\beta > -1$ and $n \in \N$,
        \begin{equation}
        \label{Jacobi.lim}
        \lim_{\alpha \to \infty} \frac{ P_n^{(\alpha, \beta)} (\cos \vartheta)}{ P_n^{(\alpha, \beta)} (1)}
         = \left( \frac{1+\cos \vartheta }{2} \right)^n, ~~~ \vartheta  \in [0, \pi].
        \end{equation}

     \item[\textup{(ii)}]  For a fixed $\beta \ge -\frac{1}{2}$ and $\vartheta \in (0, \pi]$, as $\alpha \to \infty$, the limit in (\ref{Jacobi.lim}) is uniformly for all $n \in \N$; that is,
     for any $\epsilon > 0$, there exists $A(\epsilon, \vartheta, \beta)$ such that, for any $\alpha > A(\epsilon, \vartheta, \beta)$ and $n \in \N$,
       \begin{equation}
       \label{Jacobi.lim.epsion}
       \left|    \frac{ P_n^{(\alpha, \beta)} (\cos \vartheta )}{ P_n^{(\alpha, \beta)} (1)}
                  -  \left( \frac{1+ \cos \vartheta }{2} \right)^n  \right|  < \epsilon.
       \end{equation}

        \end{itemize}

   \end{lemma}

    To prove Theorem \ref{thm3}, we need
   Lemmas 1 and 2.  With  identity (\ref{Askey1974}) taken from \cite{Askey1974},  Lemma 1 (ii) is derived from  (\ref{Askey1974}) and  Lemma 3 of \cite{MaMalyarenko2018}.
   The proof of Lemma 2 (ii) is given in Subsection \ref{lemma2.proof}, while   limit (\ref{Jacobi.lim})  in Lemma 2 (i) is from (18. 6.2) of \cite{Olver2010}.

 \section{ Isotropic covaraince matrix functions on $\Md$  generated from those in the Euclidean space  }

   For an $m \times m$ symmetric matrix function $\bC (\vartheta)$ with all entries continuous  on $[0, \infty)$,
                                      in this section we show that  it makes  
                              $\bC \left(    \rho_{\tiny{\S^{d}}}  (\x_1, \x_2)  \right) $  and 
  $\bC \left(  \rho_{\tiny{\mathbb{P}^{d}( \R)}} (\x_1, \x_2) \right)$    isotropic covariance matrix  functions on  $\S^d$ and $\mathbb{P}^d (\R)$,   respectively, 
                                      if it is compactly supported and it makes
                                      $\bC ( \| \x_1 -\x_2 \|)$ an isotropic covariance matrix function in  $\R^d$, whenever $d$ is odd.

   An $m$-variate stationary random field $\{ \bZ (\x), \x \in \R^d \}$ is said to be isotropic, if its covariance matrix function
   $\cov (\bZ (\x_1), \bZ (\x_2))$ depends only on the Euclidean distance $\| \x_1 - \x_2 \|$ between two points $\x_1$ and $\x_2$
   in $\R^d$.  When   $\{ \bZ (\x), \x \in \R^d \}$ is mean square continuous,  $\cov (\bZ (\x_1), \bZ (\x_2))$  is continuous in $\R^d$ and possesses an integral representation
   \cite{WangDuMa2014},
     $$ \cov (\bZ (\x_1), \bZ (\x_2)) =  \int_0^\infty  \Omega_d ( \| \x_1-\x_2 \| \omega) d \mathbf{F} (\omega),  ~~~~~~ \x_1, \x_2 \in \R^d,  $$
  where $\mathbf{F}( \omega), \omega \in [0, \infty ),$ is an $m \times m$  right-continuous, bounded matrix function with
  $ \mathbf{F}(0-) = \mathbf{0}$,  $\mathbf{F} ( \omega_2) - \mathbf{F} ( \omega_1)$ is positive definite for every pair of $\omega_1 $ and $\omega_2$ with $0 \le \omega_1 \le \omega_2$,
    $$   \Omega_d ( \omega ) =
                                      2^{\frac{d}{2}-1} \Gamma \left( \frac{d}{2} \right) \omega^{-\frac{d}{2}+1} J_{\frac{d}{2}-1} (\omega),
                                      ~   \omega \ge 0,  $$
                                      and $J_\nu (x)$ is the Bessel function of the first kind \cite{Szego1975}.
                                      For an integer or positive  order $\nu$,   $J_\nu(x)$ possesses 
        a series representation
  $$ J_{\nu} (x) = \left( \frac{x}{2} \right)^\nu
                   \sum_{k=0}^\infty \frac{(-1)^k}{k! \Gamma (\nu+k+1)}
                    \left( \frac{x}{2} \right)^{2 k},  ~~~~ x  > 0. $$

                                      \begin{theorem}
                                      \label{thm4}
                                      Suppose that $\bC (\vartheta)$ is an $m \times m$  symmetric matrix function on $[0, \infty)$ and all its entries are continuous on $[0, \infty)$.  For an odd $d$, if $\bC ( \| \x_1-\x_2 \|)$ is an isotropic covariance matrix function in $\R^d$, and, if all entries of $\bC (\vartheta)$ are compactly supported with
                                  $$ C_{ij} (\vartheta) = 0,   ~~~~ \vartheta    \ge \pi,  ~~~~ i, j = 1, \ldots, m, $$
                                  then    $\bC \left(  \rho_{\tiny{\S^d}}  (\x_1, \x_2) \right)$ is an isotropic covariance matrix function on $\S^d$,
                                  and   $\bC \left(  \rho_{\tiny{\mathbb{P}^d( \R)}} (\x_1, \x_2) \right)$ is an isotropic covariance matrix function on 
                                  $\mathbb{P}^d( \R)$.

                                      \end{theorem}

                                          It is not clear whether a similar result holds for an even integer $d$. Nevertheless,  the following corollary is a consequence
 of Theorem \ref{thm4}, since an isotropic covariance matrix function in $\R^d$ is also an isotropic covariance matrix function in $\R^{d-1}$ ($d \ge 2$),
 and $d-1$ is odd  for an even $d$.

\begin{corollary}
   Let  $\bC( \vartheta)$ be as in Theorem \ref{thm4}.
   For an even integer $d$, if $\bC(\| \x_1-\x_2 \|)$ is an isotropic covariance matrix function in $\R^d$, then
    $\bC \left(  \rho_{\tiny{\S^{d-1}}}  (\x_1, \x_2) \right)$ is an isotropic covariance matrix function on $\S^{d-1}$,
                                  and   $\bC \left(  \rho_{\tiny{\mathbb{P}^{d-1}( \R)}} (\x_1, \x_2) \right)$ is an isotropic covariance matrix function on 
                                  $\mathbb{P}^{d-1}( \R)$.
\end{corollary}

   The requirement that  $\bC( \vartheta)$ vanishes over $[ \pi, \infty)$ is not crucial in Theorem \ref{thm4},  since it is always possible to change the scale  for
   a compactly supported function. This results in the following corollary.

\begin{corollary}
    Suppose that  all entries of $\bC( \vartheta)$ are   continuous   on $[0, \infty)$, and
    $$C_{ij}( \vartheta)=0,   ~~~ \vartheta \ge  l,   ~ i, j  = 1, \ldots, m, $$
    where $l$ is a positive constant.
   For an odd integer $d$, if $\bC(\| \x_1-\x_2 \|)$ is an isotropic covariance matrix  function in $\R^d$, then
  $\bC \left(  \frac{l}{\pi}  \rho_{\tiny{\S^{d}}}  (\x_1, \x_2)  \right) $ is an isotropic covariance matrix  function on $\S^{d}$, and
  $\bC \left(  \rho_{\tiny{\mathbb{P}^{d}( \R)}} (\x_1, \x_2) \right)$ is an isotropic covariance matrix function on   $\mathbb{P}^d (\R)$.
 \end{corollary}

                                      Theorem \ref{thm4}, which  contains Theorems 3 and 4  of  \cite{Ma2016a} as special cases where $d=1, 3$,
                                       is conjectured in \cite{Ma2016a} with the comment that ``A difficulty arises when one deals with the connection between the two bases, the Bessel functions for $\R^d$ and ultraspherical polynomials for $\S^d$''.
                                       Such a difficulty is overcome in Theorem \ref{thm5},  where identity (\ref{thm5.eq})  builds a useful
                                      connection between an integral with respect to Jacobi polynomials and an integral with respect to the Bessel function, observing that 
  the right-hand side of (\ref{thm5.eq}) is related to the Fourier transform of the isotropic function $g( \| \x \|),  \x \in \R^d$.
 In the scalar case $m=1$,  Theorem \ref{thm4} is proved on $\S^d$ via another approach and is conjectured on $\Md$ in \cite{NieMa2019}, with an interesting example in \cite{Xu2017}.  

                                      \begin{theorem}
                                      \label{thm5}
                                      Suppose that $g(x)$ is  a continuous function on $[0, \pi]$, and that  $\alpha+\frac{1}{2}$ is a nonnegative integer.

                                      \begin{itemize}
                                      \item[\textup{(i)}]
                                      For a nonnegative integer   $\beta+\frac{1}{2}$,
                                      there is a number  $\xi_n \in [ n,  n+\alpha+\beta+1]$ such that
                                      \begin{equation}
                                      \label{thm5.eq}
                                      \begin{array}{lll}
                              &   &     \int_0^\pi g(\vartheta)   P_n^{(\alpha, \beta) } \left( \cos  \vartheta \right)
                 \sin^{2 \alpha+1}  \left( \frac{\vartheta}{2} \right) \cos^{2 \beta+1}  \left(  \frac{\vartheta}{2}  \right) d \vartheta  \\
                  &  = &   \frac{ \Gamma (n+\alpha+1)} { 2^{\alpha+1} n!} \int_0^\pi  \frac{J_\alpha (\xi_n x)}{\xi_n^\alpha} g(x) x^{\alpha+1} dx,   ~~~~~ n \in \N_0.
                     \end{array}
                                      \end{equation}

                                      \item[\textup{(ii)}]
                                      If the cosine series of  $g(x)$ converges at $x=0$, and
                                         \begin{equation}
                                         \label{thm5.ineq}
                                            \int_0^\pi  J_\alpha (\omega x) g(x) x^{\alpha+1} dx \ge 0,   ~~~~~ \omega \ge 0,
                                           \end{equation}
                                         then, for each $\beta$ with $\beta +\frac{1}{2} \in \N_0$,
                                             \begin{equation}
                                             \label{thm5.ineq2}
                                              h_n^{(\alpha, \beta)} =  \int_0^\pi g(\vartheta) P_n^{(\alpha, \beta) } \left( \cos  \vartheta \right)
                 \sin^{2 \alpha+1}  \left( \frac{\vartheta}{2} \right) \cos^{2 \beta+1}  \left(  \frac{\vartheta}{2}  \right) d \vartheta \ge 0,  ~~ n \in \N_0,
                                             \end{equation}
                                             the infinite series $\sum\limits_{n=0}^\infty  n^{\alpha+1} h_n^{(\alpha, \beta)}$ converges, and
                                             $g(x)$ can be written as the Jacobi series
                                              $$ g (x) = \sum_{n=0}^\infty  \frac{n! (2n+\alpha+\beta+1)  \Gamma (n+\alpha+\beta+1)}{ \Gamma (n+\alpha+1) \Gamma (n+\beta+1)}     h_n^{(\alpha, \beta)} P_n^{(\alpha, \beta)} (\cos x), ~~~~~ 0 \le x \le \pi. $$

                                      \end{itemize}
                                      \end{theorem}

         In a particular case where $\alpha =\beta =-\frac{1}{2}$, (\ref{thm5.eq}) holds with $\xi_n =n$.
  As a likely explanation for  why identity (\ref{thm5.eq}) works well for a positive integer $\alpha+\frac{1}{2}$,
      $\left( \frac{\pi}{2 x} \right)^{\frac{1}{2}} J_\alpha (x)$ is the spherical Bessel function of the first kind \cite{Olver2010} and is a linear combination of $\sin x$, $\cos x$, and rational functions, according to   (10.49.2) of \cite{Olver2010}.
                                  This may lead to its connection to       $P_n^{(\alpha, \beta) } \left( \cos  x \right)$, which is simply  a polynomial of $\cos x$.

 \section{Proofs}

  \subsection{Proof of Theorem \ref{thm1}}

 In the particular case $d=1$, $\Md = \S^1$, and Theorem \ref{thm1}  is known \cite{Ma2016a}, \cite{Ma2017}.
 For $d \ge 2$,
   it suffices to  verify the equivalnece bewteen (ii) and (iii), while  the equivalence between (i) and (ii) is shown in \cite{MaMalyarenko2018}.

 (ii) $\Longrightarrow$ (iii).   Suppose that $\bC (\vartheta)$ is of the form (\ref{cov.mf1}).  Making the transform $u =\cos \vartheta$, we obtain
      \begin{eqnarray*}
    \mathbf{H}_n^{(\alpha, \beta) }     &  =  & \int_0^\pi \bC (\vartheta)  P_n^{(\alpha, \beta) } \left( \cos  \vartheta \right)
                 \sin^{2 \alpha+1}  \left( \frac{\vartheta}{2} \right) \cos^{2 \beta+1}  \left(  \frac{\vartheta}{2}  \right) d \vartheta \\
   &  =  & \int_0^\pi  \left(  \sum_{k=0}^\infty  \mathbf{B}_k  P_k^{(\alpha, \beta) } \left( \cos  \vartheta \right)  \right)   P_n^{(\alpha, \beta) } \left( \cos  \vartheta \right)
                 \sin^{2 \alpha+1}  \left( \frac{\vartheta}{2} \right) \cos^{2 \beta+1}  \left(  \frac{\vartheta}{2}  \right) d \vartheta \\
       & = &  \sum_{k=0}^\infty  \mathbf{B}_k  \int_0^\pi  P_k^{(\alpha, \beta) } \left( \cos  \vartheta \right)
                  P_n^{(\alpha, \beta) } \left( \cos  \vartheta \right) \sin^{2 \alpha+1}  \left(  \frac{\vartheta}{2}  \right)
                    \cos^{2 \beta+1}  \left(  \frac{\vartheta}{2}   \right) d \vartheta \\
       & = & 2^{- (\alpha+\beta+1)}  \sum_{k=0}^\infty  \mathbf{B}_k  \int_{-1}^1  P_k^{(\alpha, \beta) } \left( u  \right)
                  P_n^{(\alpha, \beta) } \left(   u \right)  (1-u)^\alpha (1+u)^\beta  du  \\
      & = &  \frac{ \Gamma (n+\alpha+1) \Gamma (n+\beta+1)}{n!  (2n+\alpha+\beta+1) \Gamma (n+\alpha+\beta+1)}  \mathbf{B}_n,      ~~~~~ n \in \N_0,
      \end{eqnarray*}
      where the exchange between the integral and the infinite summation is ensured by the convergence of $\sum\limits_{k=0}^\infty \mathbf{B}_k P_k^{(\alpha, \beta)} (1)$, and   the last equality is due to the following orthogonal property of the Jacobi polynomials \cite{Szego1975},
     \begin{equation}
     \label{Jacobi.Orthogonal}
        \int_{-1}^1 P^{(\alpha, \beta)}_i (x) P^{(\alpha, \beta)}_j (x)  (1-x)^\alpha (1+x)^\beta  d x
                =  \left\{
                \begin{array}{ll}
                \frac{2^{\alpha+\beta+1} }{2 j +\alpha+\beta+1} \frac{\Gamma (j+\alpha+1) \Gamma (j+\beta+1)}{ j! \Gamma ( j +\alpha+\beta+1) },
                 ~   &  ~ i =j, \\
                 0, ~ & ~ i \neq j,
                 \end{array}   \right.
    \end{equation}
    for each pair of $\alpha>-1$ and $\beta>-1$.

                 The matrix $\mathbf{H}_n^{(\alpha, \beta) } $ is positive definite, since $\mathbf{B}_n$ is so.
                 The convergence of  $\sum\limits_{k=0}^\infty k^\alpha \mathbf{B}_k$ implies that that of
                 $\sum\limits_{k=0}^\infty k^{\alpha+1} \mathbf{H}_k^{(\alpha, \beta) } $, since
                 $ \lim\limits_{k \to \infty} \frac{ \Gamma (k+\alpha+1) \Gamma (k+\beta+1)}{k! \Gamma (k+\alpha+\beta+1)}  =1. $

   (iii) $\Longrightarrow$ (ii).  If $\mathbf{H}_n^{(\alpha, \beta) } $ ($n \in \N_0$) are positive definite, then so are
         $$ \mathbf{B}_n = \frac{n! (2n+\alpha+\beta+1)  \Gamma (n+\alpha+\beta+1)}{ \Gamma (n+\alpha+1) \Gamma (n+\beta+1)} \mathbf{H}_n^{(\alpha, \beta) } , ~~~ n \in \N_0. $$
      The convergence of  $\sum\limits_{n=0}^\infty n^{\alpha+1} \mathbf{H}_n^{(\alpha, \beta) } $ implies those of
                 $\sum\limits_{n=0}^\infty n^\alpha \mathbf{B}_n$,
  $\sum\limits_{n=0}^\infty \mathbf{B}_n P_n^{(\alpha, \beta)} (1)$, and the infinite series at the right-hand side of  (\ref{cov.mf2}),
  which converges to $\bC (\vartheta)$ uniformly over $[0, \pi]$.

  \subsection{Proof of Theorem \ref{thm2}}
  
  (i) For an odd $d \ge 3$, $\alpha+\frac{1}{2} = \frac{d-1}{2}$ is a positive integer.  In (\ref{Jacobi.identity2}) taking $\beta =\alpha$ and substituting $x$ by $-x$, from identity $P_n^{(\alpha, \beta)} (-x) = (-1)^n P_n^{(\beta, \alpha)} (x)$ we obtain
  $$  (1+x)^{ \frac{d-1}{2} }   P_n^{ \left( \frac{d-2}{2}, \frac{d-2}{2} \right)} (x) = \sum_{j=0}^{\frac{d-1}{2} }    a_j^{ \left(  - \frac{1}{2} \right)} (n) P_{n+j}^{ \left( \frac{d-1}{2}, -\frac{1}{2} \right)} (x), ~~~ x \in \R, ~ n \in \N_0, $$
  and, from (\ref{thm1.eq1}),
      \begin{eqnarray*}
      &  &  \mathbf{H}_n^{ \left( \frac{d-2}{2}, \frac{d-2}{2} \right)}
         =  \int_0^\pi \bC (\vartheta)  P_n^{ \left(  \frac{d-2}{2}, \frac{d-2}{2} \right)  } \left( \cos  \vartheta \right)
                 \sin^{d-1}  \left( \frac{\vartheta}{2} \right) \cos^{d-1}  \left(  \frac{\vartheta}{2}  \right) d \vartheta \\
     &  = & 2^{1-d}  \int_0^\pi \bC (\vartheta)  (1+\cos \vartheta)^{\frac{d-1}{2}} P_n^{ \left(  \frac{d-2}{2}, \frac{d-2}{2} \right)  } \left( \cos  \vartheta \right)
                 \sin^{d-1}  \left( \frac{\vartheta}{2} \right)  d \vartheta \\
     &  = & 2^{1-d}   \sum_{j=0}^{\frac{d-1}{2} }    a_j^{ \left(  - \frac{1}{2} \right)} (n)
                 \int_0^\pi \bC (\vartheta) P_{n+j}^{ \left( \frac{d-1}{2}, -\frac{1}{2} \right)} ( \cos \vartheta )
                  \sin^{d-1}  \left( \frac{\vartheta}{2} \right)  d \vartheta \\
      &  = & 2^{1-d}   \sum_{j=0}^{\frac{d-1}{2} }    a_j^{ \left(  - \frac{1}{2} \right)} (n)   \mathbf{H}_n^{ \left( \frac{d-2}{2}, -\frac{1}{2} \right)},
                  ~~~~~~ n \in \N_0,
                               \end{eqnarray*}
                               where  the positive constant $a_j^{ \left(  - \frac{1}{2} \right)} (n)$ is given by (\ref{Jacobi.identity2.coeff}) with $\alpha=\beta = \frac{d-2}{2}$.

    If $\bC \left(   \rho_{\tiny{\mathbb{P}^{d}( \R)}} (\x_1, \x_2)  \right)$ is an isotropic covariance matrix function on $\mathbb{P}^d (\R)$, then, by Theorem \ref{thm1},
    $\mathbf{H}_n^{ \left( \frac{d-2}{2}, -\frac{1}{2} \right)}$ is positive definite. So is $\mathbf{H}_n^{ \left( \frac{d-2}{2}, \frac{d-2}{2} \right)}$, $n \in \N_0$. The convergence of $\sum\limits_{n=0}^\infty  n^d \mathbf{H}_n^{ \left( \frac{d-2}{2},  -\frac{1}{2} \right) }  $
    implies that of  $\sum\limits_{n=0}^\infty  n^d \mathbf{H}_n^{ \left( \frac{d-2}{2},  \frac{d-2}{2} \right) } .$
    Thus,  $\bC \left(  \rho_{\tiny{\S^d}} (\x_1, \x_2) \right)$ is an isotropic covariance matrix function on $\S^d$, by Theorem \ref{thm1}.

    For an even $d$, if  $\bC \left(   \rho_{\tiny{\mathbb{P}^{d}( \R)}} (\x_1, \x_2)  \right)$
     is an isotropic covariance matrix function on $\mathbb{P}^d (\R)$, then it is
    an isotropic covariance matrix function on $\mathbb{P}^{d-1} (\R)$, with $d-1$ being odd, and,
    consequently,  $\bC \left(  \rho_{\tiny{\S^{d-1}}} (\x_1, \x_2) \right)$ is   an isotropic covariance matrix function on $\S^{d-1}$.

   (ii)   For an even $d \ge 4$, $\alpha = \frac{d-2}{2}$ is an even integer.  Substituting $x$ by $-x$, (\ref{Jacobi.identity1}) becomes
      \begin{equation}
      \label{thm2.proof}
        (1+x)^\alpha  P_n^{(\beta, \alpha)} (x) = \sum_{j=0}^\alpha  a_j^{(0)}  (n) P_{n+j}^{( \beta, 0)} (x), ~~~ x \in \R, ~ n \in \N_0.
       \end{equation}
       For $\beta =\alpha = \frac{d-2}{2}$,  it follows from
   (\ref{thm1.eq1}) and (\ref{thm2.proof}) that
      \begin{eqnarray*}
      &  &  \mathbf{H}_n^{ \left( \frac{d-2}{2}, \frac{d-2}{2} \right)}
         =  2^{2-d}  \int_0^\pi \bC (\vartheta)  (1+\cos \vartheta)^{\frac{d-2}{2}} P_n^{ \left(  \frac{d-2}{2}, \frac{d-2}{2} \right)  } \left( \cos  \vartheta \right)
                 \sin^{d-1}  \left( \frac{\vartheta}{2} \right)  \cos  \left( \frac{\vartheta}{2} \right)  d \vartheta \\
     &  = & 2^{2-d}   \sum_{j=0}^{\frac{d-2}{2} }    a_j^{ \left(  0 \right)} (n)
                 \int_0^\pi \bC (\vartheta) P_{n+j}^{ \left( \frac{d-1}{2}, 0 \right)} ( \cos \vartheta )
                  \sin^{d-1}  \left( \frac{\vartheta}{2} \right)  \cos  \left( \frac{\vartheta}{2} \right)  d \vartheta \\
      &  = & 2^{2-d}   \sum_{j=0}^{\frac{d-2}{2} }    a_j^{ \left(  0 \right)} (n)   \mathbf{H}_n^{ \left( \frac{d-2}{2},  0 \right)},
                  ~~~~~~ n \in \N_0.
                               \end{eqnarray*}
    If  $\bC \left(   \rho_{\tiny{\mathbb{P}^{d}( \C)}} (\x_1, \x_2)  \right)$ is an isotropic covariance matrix function on $\mathbb{P}^d (\C)$, then, by Theorem \ref{thm1},
    $\mathbf{H}_n^{ \left( \frac{d-2}{2}, 0 \right)}$ is positive definite. So is $\mathbf{H}_n^{ \left( \frac{d-2}{2}, \frac{d-2}{2} \right)}$, $n \in \N_0$.
    By Theorem \ref{thm1},  $\bC \left(  \rho_{\tiny{\S^d}} (\x_1, \x_2) \right)$ is an isotropic covariance matrix function on $\S^d$.

    For $d = 8, 12, \ldots,$ if  $\bC \left(   \rho_{\tiny{\mathbb{P}^{d}( \C)}} (\x_1, \x_2)  \right)$  is an isotropic covariance matrix function on $\mathbb{P}^d (\C)$, then
    $\mathbf{H}_n^{ \left( \frac{d-2}{2},  0 \right)}$ is positive definite, by Theorem \ref{thm1}.
    So  is $\mathbf{H}_n^{ \left( \frac{d-2}{2}, 1 \right)}$, $n \in \N_0$, by identity (\ref{H.n.idenity2}).
    As a result,  $\bC \left(   \rho_{\tiny{\mathbb{P}^{d}( \H)}} (\x_1, \x_2)  \right)$  is an isotropic covariance matrix function on $\mathbb{P}^d (\H)$.

     (iii)  It can be derived in a way similar to the proof of Part (ii).

     (iv)  Since $\bC \left(  \rho_{\tiny{\S^d}} (\x_1, \x_2) \right)$ is an isotropic covariance matrix function on $\S^d$, $\bC (\vartheta)$ is of the form (\ref{cov.mf2}) with
   $\alpha =\beta =\frac{d-2}{2}$, and, thus,
      \begin{eqnarray*}
      &  &  \bC \left( \frac{\vartheta}{2} \right)+ \bC \left( \pi- \frac{\vartheta}{2} \right)  \\
      & = &  \sum_{n=0}^\infty  \mathbf{B}_n \left\{
                 P_n^{ \left( \frac{d-2}{2}, \frac{d-2}{2} \right)}  \left( \cos \frac{\vartheta}{2} \right)
                 + P_n^{ \left( \frac{d-2}{2}, \frac{d-2}{2} \right)}  \left( -\cos \frac{\vartheta}{2} \right) \right\} \\
    & = &  \sum_{n=0}^\infty  \mathbf{B}_n \left\{
                 P_n^{ \left( \frac{d-2}{2}, \frac{d-2}{2} \right)}  \left( \cos \frac{\vartheta}{2} \right)
                 +(-1)^n P_n^{ \left( \frac{d-2}{2}, \frac{d-2}{2} \right)}  \left( \cos \frac{\vartheta}{2} \right) \right\} \\
    & = & 2  \sum_{n=0}^\infty  \mathbf{B}_{2n}
                 P_{2n}^{ \left( \frac{d-2}{2}, \frac{d-2}{2} \right)}  \left( \cos \frac{\vartheta}{2} \right)   \\
   & = & 2  \sum_{n=0}^\infty  \frac{ \Gamma \left( 2n+\frac{d}{2} \right) \Gamma (n+1) \mathbf{B}_{2n}}{
                                                       \Gamma \left( n+\frac{d}{2} \right) \Gamma (2n+1)}
                 P_{n}^{ \left( \frac{d-2}{2}, -\frac{1}{2} \right)}  \left( \cos \vartheta \right),
                                             \end{eqnarray*}
   where the second and the last  equalities follow from identities (4.1.3)
   and (4.1.5) of \cite{Szego1975}, respectively.  It follows from $\frac{\Gamma (n+\kappa+1)}{\Gamma (n+1)} \sim n^\kappa$ that
       $$  \lim_{n \to \infty} \frac{ \Gamma \left( 2n+\frac{d}{2} \right) \Gamma (n+1) }{
                                                       \Gamma \left( n+\frac{d}{2} \right) \Gamma (2n+1)} = 2^{\frac{d}{2}-1}, $$
                                                       and the convergence of $\sum\limits_{n=0}^\infty n^{\frac{d-2}{2}} \mathbf{B}_n$ implies
                                                       that of $\sum\limits_{n=0}^\infty  \frac{ \Gamma \left( 2n+\frac{d}{2} \right) \Gamma (n+1) \mathbf{B}_{2n}}{
                                                       \Gamma \left( n+\frac{d}{2} \right) \Gamma (2n+1)}
                 P_{n}^{ \left( \frac{d-2}{2}, -\frac{1}{2} \right)} (1)$.
                 By Theorem \ref{thm1},
                   $  \bC \left( \frac{\rho_{\tiny{\mathbb{P}^{d}( \R)}} (\x_1, \x_2) }{2} \right)+ \bC \left( \pi- \frac{\rho_{\tiny{\mathbb{P}^{d}( \R)}} (\x_1, \x_2)}{2} \right)  $ is an isotropic covariance matrix function on $\mathbb{P}^d (\R)$.

                   Similarly,  it follows from  identities (4.1.3)
   and (4.1.5) of \cite{Szego1975} that
                   \begin{eqnarray*}
      &  &  \bC \left( \frac{\vartheta}{2} \right)- \bC \left( \pi- \frac{\vartheta}{2} \right)  
        =  2  \sum_{n=0}^\infty  \mathbf{B}_{2n+1}
                 P_{2n+1}^{ \left( \frac{d-2}{2}, \frac{d-2}{2} \right)}  \left( \cos \frac{\vartheta}{2} \right)   \\
   & = & 2  \cos  \left( \frac{\vartheta}{2} \right)    \sum_{n=0}^\infty  \frac{ \Gamma \left( 2n+\frac{d}{2} +1\right) \Gamma (n+1) \mathbf{B}_{2n+1}}{
                                                       \Gamma \left( n+\frac{d}{2} \right) \Gamma (2n+2)}
                 P_{n}^{ \left( \frac{d-2}{2},  \frac{1}{2} \right)}  \left( \cos \vartheta \right),
                                             \end{eqnarray*}
and, from (\ref{Szego4.5.4.dual}),
      \begin{eqnarray*}
      &  &  \left( \bC \left( \frac{\vartheta}{2} \right)- \bC \left( \pi- \frac{\vartheta}{2} \right)  \right) \cos \left( \frac{\vartheta}{2} \right)  \\
       & = & \frac{1}{2}  \sum_{n=0}^\infty  \frac{ \Gamma \left( 2n+\frac{d}{2} +1\right) \Gamma (n+1) \mathbf{B}_{2n+1}}{
                                                       \Gamma \left( n+\frac{d}{2} \right) \Gamma (2n+2)}
                (1+\cos \vartheta)   P_{n}^{ \left( \frac{d-2}{2},  \frac{1}{2} \right)}  \left( \cos \vartheta \right) \\
       & = &   \sum_{n=0}^\infty  \frac{ \Gamma \left( 2n+\frac{d}{2} +1\right) \Gamma (n+1) \mathbf{B}_{2n+1}}{
                                                       \Gamma \left( n+\frac{d}{2} \right) \Gamma (2n+2) (4n +d+1)}
              \left( (2n+1)   P_{n}^{ \left( \frac{d-2}{2},  -\frac{1}{2} \right)}  \left( \cos \vartheta \right)  \right. \\
        &  &  \left.       + 2( n+1)  P_{n+1}^{ \left( \frac{d-2}{2},  -\frac{1}{2} \right)}  \left( \cos \vartheta \right) \right),
                                             \end{eqnarray*}
 which implies  that
 $ \left( \bC \left( \frac{\rho_{\tiny{\mathbb{P}^{d}( \R)}} (\x_1, \x_2)}{2} \right)- \bC \left( \pi- \frac{\rho_{\tiny{\mathbb{P}^{d}( \R)}} (\x_1, \x_2)}{2} \right)  \right) \cos \left( \frac{\rho_{\tiny{\mathbb{P}^{d}( \R)}} (\x_1, \x_2)}{2} \right)$  
 
 \noindent   is an isotropic covariance matrix function on $\mathbb{P}^d (\R)$ by Theorem \ref{thm1}.

   \subsection{Proof of Theorem \ref{thm3}}

   It suffices to establish the equivalence between  statements (i) and (ii), while  the equivalence between statements (ii) and  (iii) is due to the identity
      $ \vartheta = \frac{\pi}{2} -\arcsin ( \cos \vartheta ),$  $  \vartheta  \in [0, \pi], $ and that  between statements (ii) and  (iv) is
      made by the transform $x = \cos \frac{\vartheta}{2},  0 \le x \le 1.$

      (ii) $\Longrightarrow$ (i):  Let $\bC (\vartheta)$ take the form (\ref{thm3.eq1}).
    For each $n \in \N_0$,     $(1+\cos \rho (\x_1, \x_2) )^n$
      is an isotropic covariance function on each $\Md$, by Lemma 1 (i). So is $\bC ( \cos \rho (\x_1, \x_2))$, by Theorem 2 of
  \cite{MaMalyarenko2018}.

       (i) $\Longrightarrow$ (ii): Suppose that $\bC (\rho (\x_1, \x_2))$ is an isotropic covariance matrix function on $\mathbb{M}^\infty$.
       Then it is an isotropic covariance matrix function on each $\Md$, and, for each possible pair of $\alpha$ and $\beta$ in Table \ref{tab:1},  by Theorem 1 (ii),  $\bC (\vartheta)$ must be of the form
       \begin{equation}
       \label{thm3.proof.eq}
        \bC (\vartheta) = \sum_{n=0}^\infty \mathbf{B}_n^{(\alpha, \beta)}  \frac{P_n^{(\alpha, \beta)} (\cos \vartheta)}{ P_n^{(\alpha, \beta)} (1)},  ~~~~ \vartheta \in [ 0, \pi],
       \end{equation}
       where  $\{  \mathbf{B}_n^{(\alpha, \beta)}, n \in \N_0 \}$ is a sequence of  $m \times m $ positive definite matrices and the series
       $\sum\limits_{n=0}^\infty  \mathbf{B}_n^{(\alpha, \beta)}$ converges.

       When $\beta =\alpha = \frac{d-2}{2}$, limit  (18. 6.4) of \cite{Olver2010}  reads
       $ \lim\limits_{\alpha \to \infty} \frac{P_n^{(\alpha, \beta)} (\cos \vartheta)}{ P_n^{(\alpha, \beta)} (1)} =\cos^n \vartheta$.
       In (\ref{thm3.proof.eq})  taking $\alpha \to \infty$ and  applying Lemma 1 of \cite{Schoenberg1942} yields (see \cite{Ma2015})
            $$  \bC (\vartheta)  =  \sum_{n=0}^\infty \mathbf{B}_n \cos^n \vartheta,    ~~~~ \vartheta \in [ -\pi, \pi],  $$
            which contains (\ref{thm3.eq1}) as a special case.

      When $\beta$ is fixed as listed in Table \ref{tab:1},   we consider the scalar case $m=1$ first, under which
      (\ref{thm3.proof.eq}) reduces to
       $$ C (\vartheta) = \sum_{n=0}^\infty b_n^{(\alpha, \beta)}  \frac{P_n^{(\alpha, \beta)} (\cos \vartheta)}{ P_n^{(\alpha, \beta)} (1)},  ~~~~ \vartheta \in [ 0, \pi], $$
       where the nonnegative series $\sum\limits_{n=0}^\infty b_n^{(\alpha, \beta)}$ converges. 
       For the nonnegative convergent series $\sum\limits_{n=0}^\infty b_n^{(\alpha, \beta)}$, 
       its terms are bounded by  
         $$  0 \le  b_n^{(\alpha, \beta)} \le \sum\limits_{k=0}^\infty b_k^{(\alpha, \beta)} = C(0),   ~~~~ n \in \N_0. $$
             By Cantor't diagonal argument, there exists a subsequence $\{ \alpha_k , k \in \N \}$ and a nonnegative sequence $\{b_n, n \in \N_0 \}$ such that for any $n \in \N_0$,
\begin{equation}
\label{diag}
\lim_{k \to\infty} b^{(\alpha_k,\beta)}_n = b_n.
\end{equation}
For $\vartheta \in (0, \pi]$, we have
   \begin{eqnarray*}
    &   &   C (\vartheta) - \sum_{n=0}^\infty b_n \left( \frac{1+\cos \vartheta}{2} \right)^n  \\
    & = &  \sum_{n=0}^\infty b_n^{(\alpha_k, \beta)}  \frac{P_n^{(\alpha_i, \beta)} (\cos \vartheta)}{ P_n^{(\alpha_k, \beta)} (1)} - \sum_{n=0}^\infty b_n \left( \frac{1+\cos \vartheta}{2} \right)^n  \\
  & = &  \sum_{n=0}^\infty b_n^{(\alpha_k, \beta)} \left(  \frac{P_n^{(\alpha_k, \beta)} (\cos \vartheta)}{ P_n^{(\alpha_i, \beta)} (1)} -  \left( \frac{1+\cos \vartheta}{2} \right)^n  \right)
   +  \sum_{n=0}^\infty ( b_n^{(\alpha_k, \beta)}- b_n) \left( \frac{1+\cos \vartheta}{2} \right)^n,
   \end{eqnarray*}
   where the first  sum converges to 0 as $\alpha_k \to\infty$ by Lemma 2 (ii), and the second sum converges to 0 by the dominated convergence, since
$$ \sum_{n=0}^\infty b^{ (\alpha_k,\beta)}_n \left( \frac{1+\cos \vartheta}{2} \right)^n \le \sum_{n=0}^\infty C(0)\left(\frac{1+\cos \vartheta}{2} \right)^n
=\frac{2C(0)}{1-\cos \vartheta},$$
and  (\ref{diag}) implies
    $$ \lim_{k \to \infty} \sum_{n=0}^\infty  b_n^{(\alpha_k, \beta)} \left( \frac{1+\cos \vartheta}{2} \right)^n  = \sum_{n=0}^\infty  b_n \left( \frac{1+\cos \vartheta}{2} \right)^n. $$
    For $\vartheta=0$, (\ref{thm3.eq1}) is also valid, since its both sides  are continuous.

    In a vector case $m \ge 2$, if $\bC (\rho (\x_1, \x_2))$ is an isotropic covariance matrix function on $\mathbb{M}^\infty$, then
    $\mathbf{a}' \bC (\rho (\x_1, \x_2)) \mathbf{a}$ is an isotropic covariance  function on $\mathbb{M}^\infty$
    for  an arbitrary $\mathbf{a} \in \R^m$.
    Thus,
       $$  \mathbf{a}' \bC (\vartheta) \mathbf{a} =
           \sum_{n=0}^\infty  b_n (\mathbf{a})  (1+\cos \vartheta)^n,  ~~~~ \vartheta \in [ 0, \pi], $$
           where $\{  b_n (\mathbf{a}), n \in \N_0 \}$ is a sequence of nonnegative numbers, and $\sum\limits_{n=0}^\infty 2^n b_n (\mathbf{a})$ converges. Similarly, for  an arbitrary $\mathbf{b} \in \R^m$,
            \begin{equation}
            \label{thm3.proof.eq2}
           ( \mathbf{a}+\mathbf{b})' \bC (\vartheta) (\mathbf{a}+\mathbf{b}) =
           \sum_{n=0}^\infty  b_n (\mathbf{a}+\mathbf{b})  (1+\cos \vartheta)^n,  ~~~~ \vartheta \in [ 0, \pi],
           \end{equation}
           and
             \begin{equation}
            \label{thm3.proof.eq3}
           ( \mathbf{a}-\mathbf{b})' \bC (\vartheta) (\mathbf{a}+\mathbf{b}) =
           \sum_{n=0}^\infty  b_n (\mathbf{a}-\mathbf{b})  (1+\cos \vartheta)^n,  ~~~~ \vartheta \in [ 0, \pi].
           \end{equation}
           Taking the difference between (\ref{thm3.proof.eq2}) and  (\ref{thm3.proof.eq3}) yields
           \begin{equation}
            \label{thm3.proof.eq4}
           \mathbf{a}' \bC (\vartheta) \mathbf{b} =
           \sum_{n=0}^\infty  \frac{b_n (\mathbf{a}+\mathbf{b})- b_n (\mathbf{a}-\mathbf{b})}{2}  (1+\cos \vartheta)^n,  ~~~~ \vartheta \in [ 0, \pi],
           \end{equation}
           noticing that $\bC(\vartheta)$ is symmetric. The form (\ref{thm3.eq1}) of $\bC(\vartheta)$ and the convergence of $\sum\limits_{n=0}^\infty 2^n \mathbf{B}_n$ are obtained
           from (\ref{thm3.proof.eq4})  by taking the $i$th entry of $\mathbf{a}$
           and the $j$th entry of $\mathbf{b}$ equal to 1 and the rest being 0, for $i, j \in \{ 1, \ldots, m \}$.

     Multiplying both sides of  (\ref{thm3.eq1}) by an arbitrary $\mathbf{a} \in \R^m$ yields
    $$ \mathbf{a}'  \bC (\vartheta) \mathbf{a}
       = \sum_{n=0}^\infty  \mathbf{a}' \mathbf{B}_n \mathbf{a}   (1+\cos \vartheta )^n,  ~~~~ \vartheta \in [ 0, \pi], $$
   where the left-hand side is an isotropic covariance function on $\mathbb{M}^\infty$, so that the coefficients at the right-hand side,
    $ \mathbf{a}' \mathbf{B}_n \mathbf{a}$, have to be nonnegative; in other words, $\{   \mathbf{B}_n, n \in \N_0 \}$
    must be a sequence of positive definite matrices.

   \subsection{Proof  of Lemma 2}
   \label{lemma2.proof}

 For $\alpha>\beta>-\half$,  $\frac{ P_n^{(\alpha, \beta)} (\cos \vartheta )}{ P_n^{(\alpha, \beta)} (1)}$ admits an   integral representation
  (see,  formula (18.10.3) of \cite{Olver2010})
\begin{equation}
\label{intrep}
 \frac{ P_n^{(\alpha, \beta)} (\cos \vartheta )}{ P_n^{(\alpha, \beta)} (1)}  =
\int_0^1\int_0^\pi \left( \frac{1+\cos \vartheta}{2}-r^2 \frac{1-\cos \vartheta}{2}+ \imath r \sin \vartheta \cos \phi \right)^n h^{(\alpha,\beta)} (r,\phi)d\phi dr,
\end{equation}
where $\imath$ is the imaginary unit, and
 $$ h^{(\alpha,\beta)}(r,\phi)= \frac{(1-r^2)^{\alpha-\beta-1}r^{2\beta+1}\sin^{2\beta}\phi}{ \int_0^1\int_0^\pi (1-r^2)^{\alpha-\beta-1}r^{2\beta+1}\sin^{2\beta}\phi d \phi d r},  ~~ 0 \le r \le 1, ~  0 \le  \phi \le \pi, $$
 is a  nonnegative function with the range between 0 and 1.
 Notice that
  \begin{eqnarray*}
  &  & \left|\left(\frac{1+\cos \vartheta}{2}-r^2 \frac{1-\cos \vartheta}{2}+ \imath r \sin  \vartheta \cos\phi \right)^n
-\left(\frac{1+\cos \vartheta}{2} \right)^n \right|  \\
 &  \le &  \left(\frac{1+\cos \vartheta}{2}+r^2 \frac{1-\cos \vartheta}{2} \right)^n+\left(\frac{1+\cos \vartheta}{2} \right)^n.
 \end{eqnarray*}
 For a given $0< \vartheta \le \pi$ and $\epsilon>0$, there exists $N(\epsilon, \vartheta)$ such that for any $n>N(\epsilon, \vartheta)$ and $r<\half$,
\begin{equation}
\label{halfeps}
  \left|\left(\frac{1+\cos \vartheta}{2}-r^2 \frac{1-\cos \vartheta}{2}+ \imath r \sin  \vartheta \cos\phi \right)^n
-\left(\frac{1+\cos \vartheta}{2} \right)^n \right| < \frac{\epsilon}{2}.
\end{equation}
 On the other hand, there exists $0<\delta(\epsilon, \vartheta)<\half$ such that (\ref{halfeps}) holds for any $0\le n \le N(\epsilon, \vartheta)$ and $0 \le r<\delta(\epsilon, \vartheta)$. Therefore, (\ref{halfeps}) holds for any $0\le r<\delta(\epsilon, \vartheta)$ and $n \in \N_0$. For $\alpha>\beta>-\half$,
it follows from (\ref{intrep}) that
    \begin{eqnarray*}
     &  &  \left|  \frac{ P_n^{(\alpha, \beta)} (\cos \vartheta )}{ P_n^{(\alpha, \beta)} (1)} - \left( \frac{1+ \cos \vartheta }{2}  \right)^n  \right|     \\
     & = &  \left|   \int_0^1\int_0^\pi \left(  \frac{1+\cos \vartheta}{2}-r^2 \frac{1-\cos \vartheta}{2}+ \imath r \sin \vartheta \cos \phi \right)^n
         - \left( \frac{1+ \cos \vartheta }{2} \right)^n   h^{(\alpha,\beta)} (r,\phi)d\phi dr \right| \\
     &  \le &  \int_0^1\int_0^\pi \left|  \left( \frac{1+\cos \vartheta}{2}-r^2 \frac{1-\cos \vartheta}{2}+ \imath r \sin \vartheta \cos \phi \right)^n
         - \left( \frac{1+ \cos \vartheta }{2}  \right)^n \right|   h^{(\alpha,\beta)} (r,\phi)d\phi dr \\
     &  =  &     \int_0^{\delta(\epsilon, \vartheta)} \int_0^\pi \left|  \left( \frac{1+\cos \vartheta}{2}-r^2 \frac{1-\cos \vartheta}{2}+ \imath r \sin \vartheta \cos \phi \right)^n
         - \left( \frac{1+ \cos \vartheta }{2}  \right)^n \right|   h^{(\alpha,\beta)} (r,\phi)d\phi dr \\
     &  & +   \int_{\delta(\epsilon, \vartheta)}^1\int_0^\pi \left|  \left( \frac{1+\cos \vartheta}{2}-r^2 \frac{1-\cos \vartheta}{2}+ \imath r \sin \vartheta \cos \phi \right)^n
         - \left( \frac{1+ \cos \vartheta }{2}  \right)^n \right|   h^{(\alpha,\beta)} (r,\phi)d\phi dr \\
    & \le & \frac{\epsilon}{2}+\frac{4\Gamma(\alpha+1)(1-\delta(\epsilon, \vartheta)^2)^{\alpha-\beta-1}}{\Gamma(\alpha-\beta)\Gamma(\beta+1)} \\
    & <  &    \epsilon,
    \end{eqnarray*}
    where the last inequality holds since it possible to find $A(\epsilon, \vartheta, \beta)$ such that,
     for $\alpha > A(\epsilon, \vartheta, \beta)$,
     $$ \frac{4\Gamma(\alpha+1)(1-\delta(\epsilon, \vartheta)^2)^{\alpha-\beta-1}}{\Gamma(\alpha-\beta)\Gamma(\beta+1)} < \frac{\epsilon}{2}. $$
     Noticing that $A(\epsilon, \vartheta, \beta)$ is finite when $\beta \to -\frac{1}{2}$,
     inequality (\ref{Jacobi.lim.epsion}) also holds for $\alpha > A(\epsilon, \vartheta, \beta)$ and $\beta = -\frac{1}{2}$.

     \subsection{Proof of Theorem \ref{thm4}}

     In case $m=1$, Theorem \ref{thm4} follows directly from Theorem 5.

     For $m \ge 2$, define $g(x) = \mathbf{a}' \bC( x) \mathbf{a},  x \ge 0$, for an arbitrary $\mathbf{a} \in \R^m$.
     Since $\bC (\| \x_1-\x_2 \|)$ is an isotropic covariance matrix function in $\R^d$, $g(x)$ satisfies inequality (\ref{thm5.ineq}) by Theorems 3.1 and 3.2 of \cite{WangDuMa2014}, so that
       inequality (\ref{thm5.ineq2}) holds for each $n \in \N_0$, {\em i.e.,}  Theorem \ref{thm1} (iii) is satisfied.
       Consequently, $\bC (\rho (\x_1, \x_2))$ is an isotropic covariance matrix function on $\S^d$ or $\mathbb{P}^d (\R)$.

     \subsection{Proof of Theorem \ref{thm5}}

 (i)   For $n \in \N_0$,   write
         $$ g_n^{(\alpha, \beta)} = \frac{n!\sqrt{\pi}}{\Gamma(n+\alpha+1)}
\int_0^\pi g ( \vartheta )P_n^{(\alpha,\beta)}(\cos  \vartheta) \sin^{2\alpha+1} \left( \frac{\vartheta}{2} \right)
        \cos^{2 \beta+1} \left( \frac{\vartheta}{2} \right) d \vartheta   $$
and
   $$ g_{ (\alpha) }  (\omega) = \frac{\sqrt{\pi}}{2^{\alpha+1}}\int_0^\pi g (x)\frac{J_\alpha( \omega x)}{\omega^\alpha}x^{\alpha+1} dx,  ~~~~
           \omega \ge 0. $$
           Then (\ref{H.n.idenity}) reads
              \begin{equation}
             \label{thm5.proof1}
              (2n+\alpha+\beta+2) g_n^{(\alpha+1, \beta)}  = g_n^{(\alpha, \beta)} - g_{n+1}^{(\alpha, \beta)}, ~~~~ n \in \N_0,
               \end{equation}
               and it follows from  the identity
$\diff{}{x}\left(\frac{J_\alpha(x)}{x^\alpha}\right)
=-\frac{J_{\alpha+1}(x)}{x^\alpha}$ that
\begin{equation}
\label{Aalpha}
g_{ (\alpha+1)}(\omega)=-\inv{2 \omega}\diff{g_{(\alpha)} (\omega)}{\omega}.
\end{equation}
What needs a proof  now is the following equivalent form of identity   (\ref{thm5.eq}),
      \begin{equation}
             \label{thm5.proof2}
              g_n^{(\alpha, \beta)} = g_{ (\alpha) }  (\xi_n),   ~~~~~~ n \in \N_0.
             \end{equation}

             In a particular case where $\alpha=\beta = -\frac{1}{2}$,  (\ref{thm5.proof2}) holds with $\xi_n =n$,
           since  $P_n^{ \left( -\frac{1}{2},  -\frac{1}{2} \right) }(\cos  \vartheta) $ $= \frac{(2n)!}{2^{2n} (n!)^2} \cos (n \vartheta)$, $J_{-\frac{1}{2}} (x) = \sqrt{ \frac{ 2}{\pi x} } \cos x$, and
               $$  g_n^{ \left( -\frac{1}{2},  -\frac{1}{2} \right)} = \int_0^\pi g(\vartheta) \cos (n \vartheta)  d \vartheta = g_{ \left( -\half \right) }  (n),   ~~~~~~ n \in \N_0. $$

               Next we verify   (\ref{thm5.proof2}) for  $\alpha>-\frac{1}{2}$ and $\beta=-\frac{1}{2}$, where $\alpha+\frac{1}{2}$ is a positive integer.
Define $h(\omega) = g_{\left( -\frac{1}{2} \right)} (\sqrt{\omega})$, $\omega \ge 0$. Then
    \begin{equation}
       \label{thm5.proof3}
        g_{ (\alpha)}(\omega) =  (-1)^{\alpha+\frac{1}{2}} \frac{d^{ \alpha+\half}}{ d \omega^{\alpha+\half} }  h (\omega^2).
     \end{equation}
  By induction on $\alpha+\frac{1}{2}$ or simply on $\alpha$, we can show that
      \begin{equation}
       \label{thm5.proof4}
       g_n^{ \left( \alpha, -\half \right)} =
        (-1)^{\alpha+\half } \left( \alpha+\half \right)!
       \,  D \left[ n^2, (n+1)^2, \ldots, \left( n +\alpha+\half \right)^2 \right] h,
        \end{equation}
        where
         $$ D[y_1, \ldots, y_k] h = \sum_{j=1}^k \frac{h (y_j)}{\prod\limits_{i=1, i \neq j}^k (y_j -y_i)} $$
        is the $(k-1)$th  divided difference of $h (x)$. Indeed, this is true for $\alpha= -\half$.
        Assuming that (\ref{thm5.proof4}) is valid for an $\alpha$, then, by identity (\ref{thm5.proof1}),
           \begin{eqnarray*}
            g_n^{ \left( \alpha+1, -\half \right)}
          & = & \frac{ g_n^{ \left( \alpha, -\half \right)}- g_{n+1}^{ \left( \alpha, -\half \right)} }{
                 2n+ \alpha +\frac{3}{2} } \\
          & = &   \frac{   (-1)^{\alpha+\half } \left( \alpha+\half \right)! }{
                   2n+ \alpha +\frac{3}{2} } \left\{
                     D \left[ n^2, (n+1)^2, \ldots, \left( n +\alpha+\half \right)^2 \right] h \right. \\
           &  &   \left.
                   - D \left[ (n+1)^2, (n+2)^2, \ldots, \left( n +\alpha+\frac{3}{2} \right)^2 \right] h  \right\}  \\
      & = &   \frac{   (-1)^{\alpha+\frac{3}{2} } \left( \alpha+\frac{3}{2} \right)! }{
                   n^2- \left( n+ \alpha +\frac{3}{2} \right)^2 } \left\{
                     D \left[ n^2, (n+1)^2, \ldots, \left( n +\alpha+\half \right)^2 \right] h \right. \\
           &  &   \left.
                   - D \left[ (n+1)^2, (n+2)^2, \ldots, \left( n +\alpha+\frac{3}{2} \right)^2 \right] h  \right\}  \\
     & = &    (-1)^{\alpha+\frac{3}{2} } \left( \alpha+\frac{3}{2} \right)!
        D \left[ n^2, (n+1)^2, \ldots, \left( n +\alpha+\frac{3}{2} \right)^2 \right] h,
         \end{eqnarray*}
         {\em i.e.,} (\ref{thm5.proof4}) is valid for  $\alpha+1$.
  Applying the mean value theorem   \cite{Boor2005} to the divided difference on the right-hand side of (\ref{thm5.proof4}), $g_n^{ \left( \alpha, -\half \right)}$  can be written as
     $$ g_n^{ \left( \alpha, -\half \right)} = (-1)^{ \alpha+\half} h^{ \left( \alpha+\half \right)} (\varsigma_n), $$
     for some $\varsigma_n \in \left[ n^2, \left( n+\alpha+\half \right)^2 \right]$.
     Comparing it with (\ref{thm5.proof3}) yields
       $ g_n^{ \left( \alpha, -\half \right)} = g_{(\alpha)} (\xi_n), $
       where $\xi_n = \sqrt{\varsigma_n} \in   \left[ n,  n+\alpha+\half  \right]$.

       Lastly, we verify  (\ref{thm5.proof2}) by induction on $\beta+\half$ or $\beta$. The case of $\beta =-\half$ has been proved. Suppose that
       (\ref{thm5.proof2}) is valid for some $\beta$. By identity (\ref{H.n.idenity2}), we obtain
          \begin{eqnarray*}
             g_n^{(\alpha, \beta+1)}
          & = &  \frac{n+\beta+1}{2n+\alpha+\beta+2} g_n^{ (\alpha,\beta)}
                     +\frac{n+\alpha+1}{2n+\alpha+\beta+2} g_{n+1}^{ (\alpha,\beta)} \\
          &=&    \frac{n+\beta+1}{2n+\alpha+\beta+2} g_{ (\alpha)} (\xi_{n_1})
                   +\frac{n+\alpha+1}{2n+\alpha+\beta+2} g_{(\alpha)} (\xi_{n_2}),
          \end{eqnarray*}
where $\xi_{n_1} \in [n,n+\alpha+\beta+1]$ and $\xi_{n_2} \in [n,n+\alpha+\beta+2]$. In other words, $g_n^{(\alpha,\beta+1)}$ is an interpolation between
$ g_{(\alpha)} (\xi_{n_1})$ and $ g_{(\alpha)} (\xi_{n_2})$. Since $g_{(\alpha)} (\omega)$ is a continuous function for integrable $g(x)$, we have
$$ g_n^{(\alpha,\beta+1)} = g_{(\alpha)} (\xi_n),$$
for some $\xi_n$ between $\xi_{n_1}$ and $\xi_{n_2}$, which resides in the interval $[n, n+\alpha+\beta+2]$.

(ii)  Under assumption (\ref{thm5.ineq}),  it follows from (\ref{thm5.eq}) that  $g_n^{(\alpha,\beta)} \ge 0$. It remains to prove that $\sum\limits_{n=0}^\infty g_n^{(\alpha,\beta)}P_n^{ (\alpha,\beta)}(1)$ is bounded. For the continuous function $g(x)$, we can define the formal Jacobi series,
$$ \hat{g}^{(\alpha,\beta)}(x)=\sum_{n=0}^\infty g_n^{(\alpha,\beta)}[g] P_n^{(\alpha,\beta)}(\cos x),$$
where
\begin{eqnarray*}
  g_n^{(\alpha,\beta)}[g] & = & \frac{n! (2n+\alpha+\beta+1)\Gamma(n+\alpha+\beta+1)}{\Gamma(n+\alpha+1)\Gamma(n+\beta+1)}  \\
      &  &  \times \int_0^\pi g(x)P_n^{(\alpha,\beta)}(\cos x)\sin^{2\alpha+1} \left( \frac{x}{2} \right)  \cos^{2\beta+1} \left( \frac{x}{2} \right)  dx,
\end{eqnarray*}
and $[g]$ indicates the dependency of $g_n^{(\alpha,\beta)}$ on $g$. For $n\in \N$, define
$$S_n^{(\alpha,\beta)}[g](x)=\sum_{k=0}^{n-1} g_k^{(\alpha,\beta)}[g]P_k^{(\alpha,\beta)}(\cos x).$$
If $g(x)\in \mathbb{P}_n(\cos x)$, the space of polynomials with degree less than $n$, then $S_n^{(\alpha,\beta)}[g](x)$ $=g(x)$.
Denote by $V^{(\alpha,\beta)}$ the set of functions $h(x)$ for which $h_n^{(\alpha,\beta)}[h] \ge 0$ for all $n \in \N_0$.
As is shown in Part (i),   $ g \in V^{ \left( \alpha,-\frac{1}{2} \right)}$. What we are going to show is that
$$S_n^{(\alpha,\beta)}[g](0)=\sum_{k=0}^{n-1} g_k^{(\alpha,\beta)}[g] P_k^{(\alpha,\beta)}(1)$$
is bounded for any $\beta+\half \in\N_0$ and $n \in \N$.

First we prove the case of $\beta=-\half$. If $\alpha=-\half$, it is obvious since the Jacobi series for $\alpha=\beta=-\half$ is the cosine series. For $\alpha \ge \half$, noticing that the Jacobi functions converge to cosine functions as $n\to\infty$,
which implies by the Riemann lemma that $\lim\limits_{n\to\infty} g_n^{(\alpha,\beta)} [g]=0$, we apply (\ref{thm5.proof1})  to obtain
$V^{(\alpha,\beta)} \subseteq V^{(\alpha-1,\beta)}.$
Setting
$$\phi(x)=S_n^{(\alpha,\beta)}[g](x) \in \mathbb{P}_n(\cos x),$$
we have $ g-\phi \in V^{(\alpha,\beta)} \subseteq V^{(\alpha-1,\beta)}$, and
 \begin{eqnarray*}
  S_n^{(\alpha-1,\beta)}[g](0) & = & S_n^{(\alpha-1,\beta)}[\phi](0)+S_n^{(\alpha-1,\beta)}[g-\phi](0)  \\
                                             & \ge &  S_n^{(\alpha-1,\beta)}[\phi](0) 
                                              =  \phi(0) 
                                              =  S_n^{(\alpha,\beta)}[\phi](0) \\
                                             & = & S_n^{(\alpha,\beta)}[g](0).
  \end{eqnarray*}
As a result,
$$ S_n^{ \left( \alpha,-\half \right) }[g](0) \le S_n^{ \left(-\half,-\half \right)}[g](0).$$
Thus, the convergence of the cosine series for $g$ at $x=0$, or equivalently, the uniform boundedness of $S_n^{ \left( -\half,-\half \right)}[g](0)$, implies that $S_n^{ \left( \alpha,-\half \right)}[g](0)$ is uniformly bounded for all $n \in \N_0$.

To see that $S_n^{(\alpha,\beta)}[g](0)$ is uniformly bounded for all $\beta \ge -\half$ and $n \in \N_0$, notice that   (\ref{thm5.proof1})  implies that for the function $\phi(x)$ defined above, $g_{n-1}^{(\alpha,\beta+1)}[\phi]\ge0$, and
$$g_k^{(\alpha,\beta+1)}[g-\phi]=0$$
for all $k<n-1$. Therefore
$$ S_n^{(\alpha,\beta)}[g](0) = S_n^{(\alpha,\beta)}[\phi](0)=S_n^{(\alpha,\beta+1)}[\phi](0)
  \ge S_{n-1}^{(\alpha,\beta+1)}[\phi](0)=S_{n-1}^{(\alpha,\beta+1)}[f](0).$$
The uniform boundedness of $S_n^{ \left( \alpha,-\half \right)}[g](0)$ results in the uniform boundedness of $S_n^{(\alpha,\beta)}[g](0)$ over all $\beta\ge -\half$ and $n \in \N_0$.


\begin{thebibliography}{10}


   \bibitem{Askey1974}
   Askey, R.:  Jacobi polynomials. I. New proofs of Koornwinder's Laplace type integral representation and Bateman's bilinear sum. SIAM J. Math. Anal. {\bf 5},  119–124  (1974)

    \bibitem{Askey1976}
 Askey, R., Bingham, N. H.:
 Gaussian processes on compact symmetric spaces.
  Z. Wahrscheinlichkeitstheorie verw. Gebiete {\bf  37},  127-143 (1976)

    \bibitem{Azevedo2017}
   Azevedo, D.,  Barbosa, V. S.:
   Covering numbers of isotropic reproducing kernels on compact two-point homogeneous spaces.
   Math. Nachr.  {\bf 290},  2444--2458 (2017)

\bibitem{Bhattacharya2012}
   Bhattacharya, A., Bhattacharya, R.
   Nonparametric Inference on Manifolds. Cambridge Univ. Press, Cambridge (2012)

\bibitem{Bingham1973}
   Bingham, N. H.:
 Positive definite functions on spheres.
   Proc. Cambridge Phil. Soc.  {\bf 73},   145-156 (1973)

  \bibitem{Bochner1941}
  Bochner, S.: Hilbert distance and positive definite functions. Ann. Math. {\bf 42}, 647--656 (1941)

  \bibitem{Boor2005}
   de Boor, C.:  Divided differences. Surv. Approx. Theory {\bf 1}, 46-69  (2005)

  \bibitem{BrownDai2005}
  Brown, G.,   Dai, F.
  Approximation of smooth functions on compact two-point
              homogeneous spaces.
              J. Funct. Anal. {\bf 220}, 401-423 (2005)




  \bibitem{Cheng2016}
 Cheng, D.,    Xiao, Y.:
  Excursion  probability of Gaussian random fields on sphere.   Bernoulli {\bf 22},  1113--1130 (2016)

  \bibitem{Cohen2012}
  Cohen, S.,  Lifshits, M. A.:
  Stationary Gaussian random fields on hyperbolic spaces and on Euclidean spheres.
   ESAIM Probab. Stat.   {\bf 16},  165-221 (2012)




\bibitem{Dovidio2014}
D'Ovidio, M.:
Coordinates changed random fields on the sphere.
 J. Stat. Phys. {\bf 154}, 1153--1176 (2014)






\bibitem{Gangolli1967}
 Gangolli, R.:  Positive definite kernels on homogeneous spaces and certain stochastic processes related
to L{\' e}vy's Brownian motion of several parameters.  Ann Inst H Poincar{\' e} B  {\bf 3},  121--226 (1967)






\bibitem{Gradshteyn2007}
Gradshteyn, I. S.,  Ryzhik,  I. M.:
Tables of  Integrals, Series, and Products, 7th edtion.
Academic Press, Amsterdam (2007)





 \bibitem{Helgason2011}
  Helgason,  S.:     Integral Geometry
and Radon Transforms.   Springer,  New York (2011)







\bibitem{Leonenko2012}
 Leonenko, N.,    Sakhno, L.:    On spectral representation of tensor random fields on the sphere.
 Stoch. Anal. Appl.  {\bf 31},   167--182 (2012)


\bibitem{Leonenko2013}
Leonenko,  N.,   Shieh,  N.:    R{\' e}nyi function for multifractal random fields.
    Fractals,   {\bf  21},  1350009, 13 pp (2013)


 \bibitem{Ma2011}
  Ma, C.:
 Vector random fields with second-order moments or second-order increments.
    Stoch. Anal. Appl.  {\bf 29},   197--215 (2011)


\bibitem{Ma2015}
 Ma, C.: Isotropic covariance matrix functions on all spheres. Math. Geosci. {\bf 47},   699–717 (2015)


 \bibitem{Ma2016a}
 Ma, C.:  Stochastic representations of isotropic vector random fields on spheres.
   Stoch. Anal. Appl. {\bf 34}, 389--403 (2016)

    \bibitem{Ma2017}
 Ma, C.:  Time varying isotropic vector random fields on spheres.
   J. Theor. Prob. {\bf 30}, 1763--1785 (2017)

   \bibitem{MaMalyarenko2018}
   Ma, C.,   Malyarenko,  A.: Time-varying  isotropic vector random fields on compact two-point homogeneous spaces.
    Accepted by  J. Theor. Prob.


\bibitem{Malyarenko2013}
  Malyarenko,  A.:
   Invariant Random Fields on Spaces with a Group Action.
 Springer, New York (2013)

 \bibitem{Malyarenko1992}
  Malyarenko,  A.,   Olenko, A.:  Multidimensional covariant random fields on commutative locally compact groups.
   Ukrainian Math. J.  {\bf  44},  1384-1389 (1992)


 \bibitem{NieMa2019}
  Nie, Z.,  Ma, C.:
  Isotropic positive definite functions on spheres generated from those in Euclidean spaces.
  Proc. Amer. Math. Soc.  {\bf 147}, 3047-3056 (2019)


\bibitem{Olver2010}
Olver, F. W. J.,   Lozier, D. W.,   Boisvert, R. F.,   Clark, C. W.:
 NIST Handbook of Mathematical Functions.  Cambridge University Press, Cambridge  (2010)


\bibitem{Patrangenaru2016}
 Patrangenaru, V.,  Ellingson, L.:
Nonparametric Statistics on Manifolds and Their Applications to Object Data Analysis.
Taylor \& Francis Group, LLC, New York (2016)


 \bibitem{Schoenberg1942}
  Schoenberg, I.:   Positive definite functions on spheres.
  Duke Math. J.  {\bf 9},   96--108 (1942)

 \bibitem{Szego1975}
  Szeg\"{o}, G.
  Orthogonal Polynomials,  4th edition. Amer. Math. Soc. Colloq. Publ., vol 23. Amer. Math. Soc.,
Providence (1975)

  \bibitem{Wang1952}
  Wang,  H.-C.: Two-point homogenous spaces.
  Ann. Math.  {\bf 55},  177--191 (1952)

  \bibitem{WangDuMa2014}
  Wang, R., Du, J.,  Ma, C.:   Covariance matrix functions of isotropic vector random fields. Comm. Statist. - Theory Med. {\bf 43}, 2081-2093 (2014)

\bibitem{Xu2017}
  Xu, Y.: Positive definite functions on the unit sphere and integrals of Jacobi polynomials.
  Proc. Amer. Math. Soc. {\bf   146}, 2039-2048 (2018)

\bibitem{Yadrenko1983}
  Yadrenko,  A. M.
  Spectral  Theory of  Random Fields.
Optimization Software, New York (1983)

\bibitem{Yaglom1961}
 Yaglom,  A. M.:
Second-order homogeneous random fields. Proc. 4th Berkeley Symp. Math. Stat. Prob.  {\bf 2},   593--622 (1961)

\bibitem{Yaglom1987}
  Yaglom,  A. M.:
  Correlation  Theory of  Stationary and  Related  Random  Functions.
 vol. I.  Springer, New York (1987)

   \end{thebibliography}
\end{document}